\documentclass{article}
\usepackage{fullpage}

\usepackage[parfill]{parskip}

\usepackage[utf8]{inputenc}
\usepackage[T1]{fontenc}

\usepackage[final]{microtype}
\microtypecontext{spacing=nonfrench}

\usepackage{amsmath,amsfonts,amssymb,mathtools,amsthm}
\allowdisplaybreaks
\usepackage{stmaryrd}

\theoremstyle{definition} 

\newtheorem{lemma}{Lemma}[section]
\newtheorem{proposition}{Proposition}[section]

\usepackage{siunitx}
\usepackage{booktabs}
\usepackage{multirow,multicol,arydshln}

\usepackage{xcolor}
\usepackage{natbib}
\usepackage[hyphens]{url}
\usepackage[breaklinks]{hyperref}
\usepackage[hyphenbreaks]{breakurl}
\usepackage[noabbrev,capitalise]{cleveref}
\crefname{equation}{}{}
\colorlet{inlinkcolor}{green!50!black}
\colorlet{exlinkcolor}{red!50!black}
\hypersetup{
colorlinks = true,
allcolors = inlinkcolor,
urlcolor = exlinkcolor,
}

\usepackage{tikz,pgfplots,pgfplotstable}
\pgfplotsset{compat=1.17}
\usepgfplotslibrary{groupplots}

\newcommand{\labrel}[2]{\overset{\mathclap{\scriptscriptstyle#1}}{#2}}

\newcommand{\tsum}[1]{{\textstyle \sum_{#1}}}
\newcommand{\iin}[1]{\llbracket #1 \rrbracket}

\newcommand{\diff}{\mathrm{d}}
\newcommand{\pheq}{\hphantom{{}={}}}
\DeclareMathOperator{\per}{per}
\newcommand{\bigO}{\mathcal{O}}

\DeclareMathOperator{\tr}{tr}
\DeclareMathOperator{\intr}{int}
\DeclareMathOperator{\cl}{cl}
\DeclareMathOperator{\dom}{dom}

\newcommand{\bbR}{\mathbb{R}}
\newcommand{\bbS}{\mathbb{S}}
\newcommand{\bbH}{\mathbb{H}}
\newcommand{\Rge}{\bbR_{\geq}}
\newcommand{\Sge}{\bbS_{\succeq}}
\newcommand{\Hge}{\bbH_{\succeq}}

\newcommand{\K}{\mathcal{K}}
\newcommand{\cQ}{\mathcal{Q}}
\newcommand{\cE}{\mathcal{E}}
\newcommand{\cS}{\mathcal{S}}

\DeclareMathOperator{\logdet}{logdet}
\DeclareMathOperator{\rtdet}{rtdet}
\DeclareMathOperator{\mmd}{MMD}
\newcommand{\kmmd}{\K_{\mmd}}
\newcommand{\krtdet}{\K_{\rtdet}}
\newcommand{\klogdet}{\K_{\logdet}}

\newcommand{\tu}{\tilde{u}}
\newcommand{\tp}{\tilde{p}}
\newcommand{\tV}{\tilde{V}}
\newcommand{\bH}{\bar{H}}
\newcommand{\bZ}{\bar{Z}}
\newcommand{\hr}{\hat{r}}
\newcommand{\hxi}{\hat{\xi}}
\newcommand{\chr}{\check{r}}
\newcommand{\ba}{\bar{a}}

\DeclareMathOperator{\Diag}{Diag}

\DeclareMathOperator{\vect}{vec}

\begin{document}

\title{\Large Conic optimization with spectral functions on Euclidean Jordan algebras\thanks{
The authors thank Arkadi Nemirovski for helpful suggestions leading to barriers for root-determinant cones.
This work has been partially funded by the National Science Foundation under grant OAC-1835443 and the Office of Naval Research under grant N00014-18-1-2079.
}}
\author{Chris Coey, Lea Kapelevich, Juan Pablo Vielma}
\date{\today}
\maketitle

\begin{abstract}

Spectral functions on Euclidean Jordan algebras arise frequently in convex models.
Despite the success of primal-dual conic interior point solvers, there has been little work on enabling direct support for \emph{spectral cones}, i.e. proper nonsymmetric cones defined from epigraphs and perspectives of spectral functions.
We propose simple logarithmically homogeneous barriers for spectral cones and we derive efficient, numerically stable procedures for evaluating barrier oracles such as inverse Hessian operators.
For two useful classes of spectral cones - the \emph{root-determinant cones} and the \emph{matrix monotone derivative cones} - we show that the barriers are self-concordant, with nearly optimal parameters.
We implement these cones and oracles in our open source solver Hypatia, and we write simple, natural formulations for four applied problems.
Our computational benchmarks demonstrate that Hypatia often solves the natural formulations more efficiently than advanced solvers such as MOSEK 9 solve equivalent extended formulations written using only the cones these solvers support.

\end{abstract}


\setcounter{tocdepth}{3}
\tableofcontents

\section{Introduction}
\label{sec:intro}

A conic problem minimizes a linear function over the intersection of an affine subspace with a Cartesian product of primitive proper cones (i.e.\ irreducible, closed, convex, pointed, and full-dimensional conic sets).
Any convex problem can be represented in conic form.
An advantage of using conic form is that under mild assumptions, a conic problem has a simple and easily checkable certificate of optimality, primal infeasibility, or dual infeasibility \citep{permenter2017solving}. 

A class of functions that commonly arise in convex optimization applications are spectral functions on Euclidean Jordan algebras such as the real vectors and real symmetric or complex Hermitian matrices.
In this context, a spectral function is a real-valued symmetric function of the (real) eigenvalues.
Examples include the geometric mean (or root-determinant), the entropy (e.g. von Neumann entropy), and the trace of the inverse (e.g. the A-optimal design criterion).
Indeed, many disciplined convex programming (DCP) functions are spectral functions \citep{grant2006disciplined,grant2014cvx}.
We define \emph{spectral cones} as proper cones defined from epigraphs of convex homogeneous spectral functions or epigraphs of perspective functions of convex spectral functions.
These cones allow simple, natural conic reformulations of a wide range of convex optimization problems.
Despite this, to our knowledge there has been little prior work enabling direct support for various spectral cones in primal-dual conic solvers.

Many conic solvers (such as CSDP \citep{borchers1999csdp}, CVXOPT \citep{andersen2011interior}, ECOS \citep{domahidi2013ecos,serrano2015algorithms}, MOSEK \citep{mosek2020modeling}, SDPA \citep{yamashita2003implementation}) implement primal-dual interior point methods (PDIPMs).
Conic PDIPMs require easily computable oracles for a logarithmically homogeneous self-concordant barrier (LHSCB) function for the proper cone of the conic problem.
Complexity analysis of idealized PDIPMs shows that they converge to $\varepsilon$ tolerance in $\bigO(\sqrt{\nu} \log (1 / \varepsilon))$ iterations, where $\nu$ is the barrier parameter of the LHSCB.

Currently, most PDIPMs are specialized for the nonnegative, second order, and positive semidefinite (PSD) cones, which are cones of squares of Euclidean Jordan algebras.
However, these symmetric cones limit modeling generality and often require the construction of large extended formulations (EFs) \citep{coey2021solving}.
For many spectral cones (e.g. for the negative entropy function), equivalent EFs using only symmetric cones do not exist, and when they do, they can be impractically large.
Nonsymmetric conic PDIPMs (for example, by \citet{nesterov1996infeasible,nesterov2012towards,skajaa2015homogeneous}) can handle a much broader class of cones.
In \citet{coey2021performance}, we generalize the method by \citet{skajaa2015homogeneous}, enhance its practical performance, and test our implementation in our open-source conic solver Hypatia \citep{coey2021solving}.

A key feature of Hypatia is its generic cone interface, which allows specifying a proper cone $\K$ by implementing a small list of oracles.
Once $\K$ is defined, both $\K$ and its dual cone $\K^\ast$ may be used in any combination with other recognized cones in Hypatia to construct conic models.
The oracles to implement are: an initial interior point $t \in \intr(\K)$, a feasibility test for the cone interior $\intr(\K)$ (and optionally for the dual cone interior $\intr(\K^\ast)$), and several derivative oracles for an LHSCB for the cone.
The LHSCB oracles needed for ideal performance are the gradient, the Hessian operator (i.e. the second order directional derivative applied once to a given direction), the inverse of the Hessian operator, and the third order directional derivative (applied twice to a given direction).
Fast and numerically stable procedures for evaluating oracles are crucial for practical performance in conic PDIPM solvers such as Hypatia.

Our first main contribution is to define simple logarithmically homogeneous barriers for spectral cones and derive efficient and numerically stable barrier oracle procedures.
For example, for the case where the spectral function of the cone is separable, we show how to apply the inverse Hessian operator of the barrier function very cheaply using a closed-form formula, without the need to compute or factorize an explicit Hessian matrix (which can be expensive and prone to numerical issues).
Similarly, for the negative log-determinant and root-determinant spectral cones, we derive highly-efficient specialized oracle procedures.

Our second main contribution is to show that for two important subclasses of spectral cones - the root-determinant cones and the matrix monotone derivative (MMD) cones - the barriers we propose are LHSCBs.
These LHSCBs have parameters that are only a small additive increment of one larger than the parameter of the LHSCB for the cone of squares domain of the cone, hence the parameters are near-optimal.
MMD cones allow modeling epigraphs of a variety of useful separable spectral functions, in particular the trace (or sum) of the negative logarithm (i.e. the negative log-determinant), negative entropy, and certain power functions.
Furthermore, the dual cones of the MMD cones allow modeling epigraphs of even more separable spectral functions, such as the trace of the inverse, exponential, and more power functions.

These two contributions enable efficient and numerically stable implementations of the MMD cone and the log-determinant and root-determinant cones in nonsymmetric conic PDIPMs.
We define these cones through Hypatia's cone interface.
Our MMD cone implementation is parametrized by both a Jordan algebra domain and an MMD function, allowing the user to define new domains and MMD functions.
An MMD function is easily specified by implementing a small set of oracles for its univariate form: the function itself, its first three derivatives, and its convex conjugate, as well as an interior point for the corresponding MMD cone.
We predefine five common MMD functions and three typical Jordan algebra domains: the real vectors, real symmetric matrices, and complex Hermitian matrices.

We use these new spectral cones in Hypatia to formulate example problems from distribution estimation, experiment design, quantum information science, and polynomial optimization.
The natural formulations (NFs) using these cones are simpler and smaller than equivalent EFs written in terms of the handful of standard cones recognized by either ECOS or MOSEK 9 (i.e. the common symmetric cones and the three-dimensional exponential and power cones).
Our computational experiments demonstrate that, across a wide range of sizes and spectral functions, Hypatia can solve the NFs faster than Hypatia, MOSEK, or ECOS can solve the equivalent EFs.
Furthermore, to illustrate the practical impact of our efficient oracle procedures, we show that our closed-form formula for the MMD cone inverse Hessian product is faster and more numerically reliable than a naive direct solve using an explicit Hessian factorization.

\subsection{Overview}
\label{sec:overview}

We describe relevant aspects of Euclidean Jordan algebras, cones of squares, and spectral decompositions in \cref{sec:jordan}.
In \cref{sec:spectral}, we define spectral functions on Euclidean Jordan algebras and give expressions for gradients and second and third order directional derivatives of spectral functions. 
We also specialize these formulae for separable spectral functions and the log-determinant case.

In \cref{sec:cones}, we define spectral function cones (and their dual cones) from epigraphs of homogenized convex spectral functions on cones of squares.
We propose simple logarithmically homogeneous barriers for these cones and describe the additional properties that must be satisfied by an LHSCB. 
We also define several barrier oracles needed by Hypatia's PDIPM.
Then in \cref{sec:epiper}, we describe fast and numerically stable procedures for these barrier oracles, using the derivative results from \cref{sec:spectral}.
We specialize the oracle procedures for cones defined from separable spectral functions and the log-determinant function.

In \cref{sec:matmono}, we define the MMD cone and its dual cone, and we give useful examples of MMD functions.
We show that for the MMD cone, our barrier function is an LHSCB.
In \cref{sec:rtdet}, we define the root-determinant cone and its dual cone, prove that our barrier is an LHSCB, and derive efficient oracle procedures.
Finally, in \cref{sec:testing}, we describe a series of applied examples over the new root-determinant, log-determinant, and MMD cones and their dual cones.
We perform computational testing to demonstrate the advantages of solving these NFs with Hypatia and to exemplify the impact of efficient oracle procedures.

\subsection{Notation}
\label{sec:notation}

We often write equation references above relation symbols to indicate the use of earlier results.
For a natural number $d$, we define the index set $\iin{d} = \{1, 2, \ldots, d\}$.
For a set $\mathcal{C}$, $\cl(\mathcal{C})$ and $\intr(\mathcal{C})$ denote the closure and interior.
$\bbR$ is the space of reals, $\Rge$ is the cone of nonnegative reals, $\bbR_{>} = \intr(\Rge)$ is the positive reals, and $\bbR^d$ is the $d$-dimensional real vectors.
$\bbS^d$ and $\bbH^d$ are the real symmetric and complex Hermitian matrices with side dimension $d$, and $\Sge^d$ and $\Hge^d$ are the corresponding cones of positive semidefinite matrices.

The inner product of $w$ and $r$ is $\langle w, r \rangle$.
The $k$th derivative of a function $f$ evaluated at a point $w$ is $\nabla^k f(w)$, which may be interpreted as an operator.
For example, the second directional derivative of $f$ at $w$ in the direction $r, r$ is $\nabla^2 f(w)[r, r] = \langle \nabla^2 f(w)[r], r \rangle$.
Often we omit the point at which the derivative is evaluated if this is clear from context.
We use subscripts for partial derivatives, for example $\nabla_w f$.

\section{Jordan algebras}
\label{sec:jordan}

Jordan algebraic concepts provide a useful and straightforward abstraction for spectral functions, cones, and our barrier results in later sections.
We follow the notation of \citet[Chapter 2]{faraut1998analysis} where possible.

An algebra over the real or complex numbers is a vector space $V$ equipped with a bilinear product $\circ: V \times V \to V$.
For $w \in V$, $w^2 \coloneqq w \circ w$.
We refer to $V$ as a \emph{Jordan algebra} if for all $w_a, w_b \in V$:
\begin{subequations}
\begin{align}
w_a \circ w_b &= w_b \circ w_a,
\\
w_a \circ (w_a^2 \circ w_b) &= w_a^2 \circ (w_a \circ w_b).
\end{align}
\label{eq:jordanaxioms}
\end{subequations}
For example, for $V = \bbR^d$, we can define $\circ$ as an elementwise multiplication, or for $V = \bbS^d$ and $V = \bbH^d$, we can let $w_a \circ w_b = \frac{1}{2} (w_a w_b + w_b w_a)$.

Given $w_a \in V$, we define the linear map $L(w_a): V \to V$ satisfying:
\begin{align}
L(w_a) w_b &= w_a \circ w_b & \forall w_b \in V.
\end{align}
Given $w \in V$, we define the linear map $P(w): V \to V$ satisfying:
\begin{align}
P(w) = 2 L(w)^2 - L(w^2).
\end{align}
$P$ is called the \emph{quadratic representation} of $V$.
In general, $P(w) \neq L(w)^2 \neq L(w^2)$ because $\circ$ need not be associative.
For example, for $V = \bbS^d$, we have $L(w_a^2) w_b = \frac{1}{2}(w_a^2 w_b + w_b w_a^2)$, $L(w_a)^2 w_b = \frac{1}{2} ( L(w_a^2) w_b + w_a w_b w_a )$, and $P(w_a) w_b = w_a w_b w_a$.

For any positive integer $k$, we have \citep[Corollary 2.3.9]{vieira2007jordan}:
\begin{align}
P(w)^k &= P(w^k).
\label{eq:Pk}
\end{align}
It is standard to assume the existence of a multiplicative identity $e$.
Note that $P(w) e = w$.
A point $w \in V$ is \emph{invertible} if and only if $L(w)$ is invertible, and the inverse of $w$ is the element $w^{-1} \in V$ such that $w^{-1} = L(w)^{-1} e$ \citep[Proposition II.2.2]{faraut1998analysis}.
\Cref{eq:Pk} also holds for $k = -1$ if $w$ is invertible \citep[Proposition II.3.1]{faraut1998analysis}.

Henceforth we consider only the finite dimensional \emph{Euclidean} Jordan algebras.
A Jordan algebra $V$ is Euclidean if $\langle w_a \circ w_b, w_c \rangle = \langle w_b, w_a \circ w_c \rangle$ for all $w_a, w_b, w_c \in V$. 

We call $\cQ$ a \emph{cone of squares} on $V$ if $\cQ = \{ w \circ w : w \in V \}$.
The cone $\cQ$ is proper (closed, convex, pointed, and solid) because $V$ is Euclidean (and therefore formally real); see \citet[Theorem 3.3]{papp2013semidefinite} and \citet[Section III.1 and Proposition VIII.4.2]{faraut1998analysis}.
In addition, $\cQ$ is \emph{self-dual} and \emph{homogeneous}; see \citet[Proposition 2.5.8]{vieira2007jordan} and \citet[Theorem III.2.1]{faraut1998analysis}.
For example, for $V = \bbS^d$, the cone of squares is $\cQ = \Sge^d$.

For convenience, we often write $a \succeq b$ instead of $a - b \in \cQ$, or $a \succ b$ instead of $a - b \in \intr(\cQ)$, where $\cQ$ is clear from context.
If $w \succ 0$, then $w$ is invertible \citep[Theorem III.2.1]{faraut1998analysis}.
Furthermore $w \succ 0$ implies that $w^{1/2}$ is well-defined and invertible, and $P(w^{1/2}) = P(w)^{1/2}$ \citet[Proposition 2.5.11]{vieira2007jordan}.
This also implies by \cref{eq:Pk} (with $k = -1$) that $P(w^{-1/2}) = P(w)^{-1/2}$.

\subsection{Spectral decomposition}
\label{sec:jordan:specdecomp}

In a Euclidean Jordan algebra $V$, an \emph{idempotent} is an element $c \in V$ such that $c^2 = c$.
Two idempotents $c_1, c_2$ are \emph{orthogonal} if $c_1 \circ c_2 = 0$.
Let $d$ be the rank of $V$.
$c_1, \ldots, c_d$ is a \emph{complete} system of orthogonal idempotents if $c_1, \ldots, c_d$ are all idempotents, pairwise orthogonal, and $\sum_{i \in \iin{d}} c_i = e$.
An idempotent is \emph{primitive} if it is non-zero and cannot be written as the sum of two orthogonal non-zero idempotents.
A \emph{Jordan frame} is a complete system of orthogonal idempotents, where each idempotent is primitive.
The number of elements in any Jordan frame is called the \emph{rank} of $V$.
For example, the rank of $\bbR^d$, $\bbS^d$, or $\bbH^d$ is $d$.

For any $w \in V$, there exist unique real numbers (not necessarily distinct) $w_1, \ldots, w_d$ and a unique Jordan frame $c_1, \ldots, c_d$ such that $w$ has the \emph{spectral decomposition} \citep[Theorem III.1.2]{faraut1998analysis}:
\begin{equation}
w = \sum_{i \in \iin{d}} w_i c_i.
\label{eq:specdecomp}
\end{equation} 
We call $w_1, \ldots, w_d$ the \emph{eigenvalues} of $w$.
The \emph{determinant} is $\det(w) = \prod_{i \in \iin{d}} w_i$ and the \emph{trace} is $\tr(w) = \sum_{i \in \iin{d}} w_i$ \citep[Section II.2, page 29]{faraut1998analysis}.
For example, for $V = \bbR^d$, the Jordan frame is the standard unit vectors and $w$ is its own vector of eigenvalues.
For $V = \bbS^d$, we can think of the Jordan frame as the rank one PSD matrices from a full symmetric eigendecomposition.

Henceforth, we define the inner product on $V$ as $\langle w_a, w_b \rangle = \tr(w_a \circ w_b)$.
Under this inner product, $P(w)$ is self-adjoint \citep[Page 27]{vieira2007jordan}.
Thus for $w \in \intr(\cQ)$ and $r_1, r_2 \in V$, we have:
\begin{equation}
\langle P(w) r_1, r_2 \rangle 
= \langle P(w^{1/2}) r_1, P(w^{1/2}) r_2 \rangle 
= \langle r_1, P(w) r_2 \rangle.
\label{eq:moveP}
\end{equation}

\subsection{Peirce decomposition}
\label{sec:jordan:peircedecomp}

We let $c_1, \ldots, c_d$ be a Jordan frame for $V$, and define for $i,j \in \iin{d}$:
\begin{subequations}
\begin{align}
V(c_i, \lambda) & \coloneqq \{ w : c_i \circ w = \lambda w \},
\\
V_{i,i} & \coloneqq V(c_i, 1) = \{ t c_i : t \in \bbR \}, 
\\
V_{i,j} & \coloneqq V(c_i, \tfrac{1}{2}) \cap V(c_j, \tfrac{1}{2}).
\end{align}
\end{subequations}
$V$ has the direct sum decomposition $V = \oplus_{i, j \in \iin{d} : i \leq j} V_{i,j}$ \citep[Theorem IV.1.3]{faraut1998analysis}.
For example, for $V = \bbS^d$, let $E_{i,j}$ be a matrix of zeros except in the $(i,j)$th position, and let $c_i = E_{i,i}$; then $V_{i,i} = \{ t E_{i,i} : t \in \bbR \}$ and $V_{i,j} = \{ t (E_{i,j} + E_{j,i}) : t \in \bbR \}$.

The \emph{Peirce decomposition} allows us to write any $r \in V$ as:
\begin{equation}
r = \sum_{i, j \in \iin{d} : i \leq j} r_{i,j} 
= \sum_{i, j \in \iin{d} : i < j} r_{i,j} + \sum_{i \in \iin{d}} r_i c_i,
\label{eq:peircedecomp}
\end{equation}
where $r_i = \langle r, c_i \rangle \in \bbR$ and $r_{i,j} \in V_{i,j}$ for all $i, j \in \iin{d}$.
Each $r_{i,j}$ is a projection of $r$ onto $V_{i,j}$, where:
\begin{subequations}
\begin{align}
r_{i,i} &=
r_i c_i
= P(c_i) r
& \forall i \in \iin{d},
\\
r_{i,j} 
&= 4 L(c_i) L(c_j) r
= 4 c_i \circ (c_j \circ r)
& \forall i,j \in \iin{d}: j \neq i.
\end{align}
\label{eq:Vijproj}
\end{subequations}
Note that $r_{i,j} = r_{j,i}$, since $L(c_i)$ and $L(c_j)$ commute \citep[Lemma IV.1.3]{faraut1998analysis}.
For example, let $c_1, \ldots, c_d$ be a Jordan frame for $V = \bbS^d$ and let $r \in V$; then $r_i = c_i r c_i$ and $r_{i,j} = c_i r c_j + c_j r c_i$, for $i, j \in \iin{d}$.

We list some facts relating to compositions of projection operators (see \citet[Theorem IV.2.2]{faraut1998analysis} and \citet[page 430]{sun2008lowner}):
\begin{subequations}
\begin{align}
L(c_i) L(c_j) L(c_k) L(c_l) &= 0 
& \forall i, j, k, l \in \iin{d} : i \neq j, k \neq l, (i,j) \neq (k,l), 
\\
L(c_i) L(c_j) P(c_k) = P(c_k) L(c_i) L(c_j) &= 0 
& \forall i, j, k \in \iin{d} : i \neq j,
\\
(4 L(c_i) L(c_j))^2 &= 4 L(c_i) L(c_j) 
& \forall i, j \in \iin{d},
\\
P(c_i)^2 &= P(c_i) 
& \forall i \in \iin{d},
\\
\sum_{i,j \in \iin{d} : i < j} 4 L(c_i) L(c_j) + \sum_{i \in \iin{d}} P(c_i) 
&= L(e).
\end{align}
\label{eq:projfacts}
\end{subequations}
Given $\lambda_{i,j} \neq 0 \in \bbR$ for $i, j \in \iin{d}$, consider an operator $\Lambda: V \to V$ of the form:
\begin{align}
\Lambda 
&\coloneqq \sum_{i,j \in \iin{d} : i < j} 4 \lambda_{i,j} L(c_i) L(c_j) + 
\sum_{i \in \iin{d}} \lambda_{i,i} P(c_i).
\end{align}
The inverse operator $\Lambda^{-1}$ is given by:
\begin{align}
\Lambda^{-1} 
= \sum_{i,j \in \iin{d} : i < j} 4 \lambda_{i,j}^{-1} L(c_i) L(c_j) + 
\sum_{i \in \iin{d}} \lambda_{i,i}^{-1} P(c_i).
\label{eq:opinv}
\end{align}
It can be verified using \cref{eq:projfacts} that for any $r \in V$, $\Lambda \Lambda^{-1} r = \Lambda^{-1} \Lambda r = r$.
For example, let $w = \sum_{i \in \iin{d}} w_i c_i \in V$ be invertible and suppose that $\lambda_{i,j} = w_i w_j$ for $i, j \in \iin{d}$; then:
\begin{subequations}
\begin{align}
\Lambda 
&= \sum_{i,j \in \iin{d} : i < j} 4 w_{i} w_j L(c_i) L(c_j) 
+ \sum_{i \in \iin{d}} w_{i}^2 P(c_i) 
= P(w),
\label{eq:Pexample}
\\
\Lambda^{-1} 
&= \sum_{i,j \in \iin{d} : i < j} 4 w_{i}^{-1} w_j^{-1} L(c_i) L(c_j) 
+ \sum_{i \in \iin{d}} w_{i}^{-2} P(c_i) 
= P(w^{-1}).
\label{eq:Lambdaex}
\end{align}
\end{subequations}

\section{Spectral functions and derivatives}
\label{sec:spectral}

Let $V$ be a Jordan algebra of rank $d$.
A real-valued function $f : \bbR^d \to \bbR$ is \emph{symmetric} if it is invariant to the order of its inputs.
A symmetric function $f$ composed with an eigenvalue map $\lambda : V \to \bbR^d$ induces a \emph{spectral function} $\varphi: V \to \bbR$ such that $\varphi(w) = f(\lambda(w))$, where $\lambda(w) = (w_1, \ldots, w_d)$ is the eigenvalue vector of $w$ \citep[Definition 8]{baes2007convexity}.
Note that $\varphi$ is convex if and only if $f$ is convex \citep{davis1957all}.

In this section, we give expressions for certain derivatives and directional derivatives of $\varphi$ that are useful for the barrier oracles we derive in \cref{sec:epiper}.
We express these derivatives at a point $w \in V$ (satisfying certain assumptions as necessary below) with spectral decomposition \cref{eq:specdecomp}, and we let the direction be $r \in V$ with Peirce decomposition \cref{eq:peircedecomp}.
The gradient is $\nabla \varphi(w) \in V$ and the second and third order directional derivatives are $\nabla^2 \varphi(w) [r] \in V$ and $\nabla^3 \varphi(w) [r, r] \in V$.
We begin with the general nonseparable case in \cref{sec:spectral:gen} before specializing for separable spectral functions in \cref{sec:spectral:sep} and finally for the important case of the negative log-determinant function in \cref{sec:spectral:log}.

\subsection{The nonseparable case}
\label{sec:spectral:gen}

Let $\nabla f$, $\nabla^2 f$, and $\nabla^3 f$ denote the derivatives of $f$ evaluated at $\lambda(w)$.
We use subindices to denote particular components of these derivatives.
According to \citet[Theorem 38]{baes2007convexity} and \citet[Theorem 4.1]{sun2008lowner}, the gradient of $\varphi$ at $w$ is:
\begin{equation}
\nabla \varphi(w) = \sum_{i \in \iin{d}} (\nabla f)_i c_i.
\label{eq:nonsepnabla1}
\end{equation}
Henceforth we assume the eigenvalues of $w$ are all distinct for simplicity. 
The second order directional derivative of $\varphi$ in direction $r$ is \citep[Theorem 4.2]{sun2008lowner}:
\begin{equation}
\nabla^2 \varphi(w) [r] = 
\sum_{i, j \in \iin{d} : i < j} \frac{(\nabla f)_i - (\nabla f)_j}{w_i - w_j} r_{i,j} + 
\sum_{i,j \in \iin{d}} (\nabla^2 f)_{i,j} r_i c_j.
\label{eq:nonsepnabla2}
\end{equation}
\Citet[Theorem 4.2]{sun2008lowner} also generalize this expression to allow for non-distinct eigenvalues.

To derive an expression for the third order directional derivative $\nabla^3 \varphi(w)[r, r]$, we let:
\begin{equation}
w(t) \coloneqq w + t r = \sum_{i \in \iin{d}} w_i(t) c_i(t),
\end{equation}
where $w_i(t)$ is the $i$th eigenvalue of $w(t)$.
Note that $\nabla^2 \varphi(w) [r] = \frac{\diff}{\diff t} \nabla \varphi(w(t)) \vert_{t=0}$ and $\nabla^3 \varphi(w) [r, r] = \frac{\diff^2}{\diff t^2} \nabla \varphi(w(t)) \vert_{t=0}$.
We let $\nabla f(t)$, $\nabla^2 f(t)$, and $\nabla^3 f(t)$ denote the derivatives of $f$ evaluated at $\lambda(w(t))$.
Due to the chain rule and \cref{eq:nonsepnabla2}:
\begin{subequations}
\begin{align}
\frac{\diff}{\diff t} \nabla \varphi (w(t))
&= \nabla^2 \varphi(w(t)) [r]
\\
&= \sum_{i,j \in \iin{d}: i < j} 
\frac{(\nabla f(t))_i - (\nabla f(t))_j}{w_i(t) - w_j(t)} r_{i,j}(t) + 
\sum_{i,j \in \iin{d}} (\nabla^2 f(t))_{i,j} r_i(t) c_j(t).
\end{align}
\label{eq:dphidt}
\end{subequations}
We differentiate \cref{eq:dphidt} once more.
From \citet[Corollary 1 and Theorem 3.3]{vieira2016derivatives} and \citet[Equation 37]{sun2008lowner}, for $i \in \iin{d}$ we have:
\begin{subequations}
\begin{align}
\frac{\diff}{\diff t} w_i(t)
&= r_i(t),
\label{eq:dldt}
\\
\frac{\diff}{\diff t} c_i(t)
&= s_i(t)
\coloneqq \sum_{j \in \iin{d}: j \neq i} \frac{r_{i,j}(t)}{w_i(t) - w_j(t)}.
\end{align}
\end{subequations}
Using the chain rule, \cref{eq:dldt} implies:
\begin{subequations}
\begin{align}
\frac{\diff}{\diff t}(\nabla f(t))_i 
&= \sum_{k \in \iin{d}} (\nabla^2 f(t))_{i,k} r_k(t)
& \forall i \in \iin{d},
\\
\frac{\diff}{\diff t}(\nabla^2 f(t))_{i,j} 
&= \sum_{k \in \iin{d}} (\nabla^3 f(t))_{i,j,k} r_k(t)
& \forall i, j \in \iin{d}.
\end{align}
\end{subequations}
Applying the chain and product rules, we have:
\begin{subequations}
\begin{align}
\frac{\diff}{\diff t} \frac{1}{w_i(t) - w_j(t)} 
&= \frac{r_j(t) - r_i(t)}{(w_i(t) - w_j(t))^2}
& \forall i,j \in \iin{d} : i \neq j,
\label{eq:rik}
\\
\frac{\diff}{\diff t} r_{i,j}(t)
&= \frac{\diff}{\diff t} ( 4 c_i(t) \circ (c_j(t) \circ r) )
\\
&= 4 c_i(t) \circ (s_j(t) \circ r) + 4 s_i(t) \circ (c_j(t) \circ r)
& \forall i,j \in \iin{d} : i \neq j,
\\
\frac{\diff}{\diff t} \langle c_i(t), r \rangle c_j(t)
&= \langle s_i(t), r \rangle c_j(t) + r_i(t) s_j(t)
& \forall i,j \in \iin{d}.
\label{eq:partials}
\end{align}
\end{subequations}
Finally, letting $s_i \coloneqq s_i(0)$ for all $i \in \iin{d}$, these results imply that:
\begin{subequations}
\begin{align}
\nabla^3 \varphi[r, r]
&= \frac{\diff^2}{\diff t^2} \nabla \varphi (w(t)) \big\vert_{t=0}
\\
\begin{split}
&= \sum_{i, j \in \iin{d}: i < j} \frac{(\nabla f)_i - (\nabla f)_j}{w_i - w_j}
\biggl( 
4 c_i \circ ( s_j \circ r) + 4 s_i \circ (c_j \circ r) - 
\frac{r_i - r_j}{w_i - w_j} r_{i,j} \biggr) + {}
\\
&\pheq \sum_{i,j,k \in \iin{d} : i < j} \frac{(\nabla^2 f)_{i,k} - 
(\nabla^2 f)_{j,k}}{w_i - w_j} r_k r_{i,j} + {}
\\
&\pheq \sum_{i, j \in \iin{d}} (\nabla^2 f)_{i,j} 
( \langle s_i, r \rangle c_j + r_i s_j ) + 
\sum_{i, j, k \in \iin{d}} (\nabla^3 f)_{i,j,k} r_i r_k c_j
\end{split}
\\
\begin{split}
&= \sum_{i, j \in \iin{d}: i < j} \frac{(\nabla f)_i - (\nabla f)_j}{w_i - w_j}
\biggl( 
4 c_i \circ ( s_j \circ r) + 4 s_i \circ (c_j \circ r) - 
\frac{r_i - r_j}{w_i - w_j} r_{i,j} \biggr) + {}
\\
&\pheq \sum_{i, j \in \iin{d}} (\nabla^2 f)_{i,j} 
(2 r_j s_i + \langle s_i, r \rangle c_j) + 
\sum_{i, j, k \in \iin{d}} (\nabla^3 f)_{i,j,k} r_i r_k c_j 
.
\end{split}
\label{eq:nonsepnabla3}
\end{align}
\end{subequations}

The derivative expressions simplify significantly for $V = \bbR^d$.
For $V = \bbS^d$, the form of \cref{eq:nonsepnabla2} is well-known \citep{faybusovich2021long} and the form of \cref{eq:nonsepnabla3} appears in \citet{sendov2007higher}.

\subsection{The separable case}
\label{sec:spectral:sep}

The spectral function $\varphi$ induced by $f$ is \emph{separable} if $f$ is a separable function, i.e. $f(\lambda) = \sum_{i \in \iin{d}} h(\lambda_i)$ for $\lambda \in \bbR^d$ and some function $h : \bbR \to \bbR$.
For convenience, if $w \in V$, we also define $h : V \to V$ as $h(w) \coloneqq \sum_{i \in \iin{d}} h(w_i) c_i$.
This allows us to write $\varphi(w) = \tr(h(w)) = \sum_{i \in \iin{d}} h(\lambda_i)$.
Note that $\varphi$ is convex if and only if $h$ is convex.
For example, if $h(w) = -\log(w)$ then $\varphi(w) = \tr(-\log(w)) = -\logdet(w)$; we consider this special case in \cref{sec:spectral:log}.

We specialize the derivatives from \cref{eq:nonsepnabla1,eq:nonsepnabla2,eq:nonsepnabla3}, maintaining the simplifying assumption of distinct eigenvalues.
Since $(\nabla^2 f)_{i,j} = (\nabla^3 f)_{i,j,k} = 0$ unless $i=j=k$, we have:
\begin{subequations}
\begin{align}
\nabla \varphi(w) 
&= \sum_{i \in \iin{d}} \nabla h(w_i) c_i,
\label{eq:n1phi}
\\
\nabla^2 \varphi(w)[r] 
&= \sum_{i,j \in \iin{d} : i < j} 
\frac{\nabla h(w_i) - \nabla h(w_j)}{w_i - w_j} 4 c_i \circ (c_j \circ r) + 
\sum_{i \in \iin{d}} \nabla^2 h(w_i) P(c_i) r,
\label{eq:n2phi}
\\
\nabla^3 \varphi(w)[r, r] 
&= \sum_{i,j \in \iin{d}: i < j} \frac{\nabla h(w_i) - \nabla h(w_j)}{w_i - w_j}
\biggl( 4 c_i \circ ( s_j \circ r) + 
4 s_i \circ (c_j \circ r) -
\frac{r_i - r_j}{w_i - w_j} r_{i,j} \biggr) + {}
\\
&\pheq \sum_{i \in \iin{d}} \nabla^2 h(w_i) 
(2 r_i s_i + \langle s_i, r \rangle c_i) +
\sum_{i \in \iin{d}} \nabla^3 h(w_i) r_i^2 c_i.
\end{align}
\label{eq:sepDphi}
\end{subequations}

\subsection{The negative log-determinant case}
\label{sec:spectral:log}

The negative log-determinant function $\varphi(w) = -\logdet(w)$ is a separable spectral function.
We let $w \succ 0$ and drop the assumption of distinct eigenvalues.
For convenience, we let $\hr \coloneqq P(w^{-1/2}) r \in V$.
First, note that (similar to \citet[Lemma 3.3.4]{vieira2007jordan}):
\begin{align}
\langle w^{-1}, r \rangle
= \langle P(w^{-1/2}) e, r \rangle
\labrel{\ref{eq:moveP}}{=} \langle e, P(w^{-1/2}) r \rangle
= \tr(\hr).
\label{eq:trhr}
\end{align}
Due to \citet[Proposition II.2.3]{faraut1998analysis}:
\begin{align}
\nabla_w (\tr(\hr)) 
= \nabla_w (w^{-1})[r] 
= -P(w^{-1}) r.
\label{eq:Dinv}
\end{align}
Adapting the result in \citet[Lemma 3.4]{faybusovich2017matrix}:
\begin{align}
\nabla_w (P(w^{-1}) r) [r]
= -2 P(w^{-1/2}) \hr^2.
\label{eq:Pinv}
\end{align}
Now, the gradient of $\varphi$ is \citep[Propositions III.4.2(ii)]{faraut1998analysis}:
\begin{align}
\nabla \varphi(w) & = -w^{-1},
\label{eq:Dpsi}
\end{align}
so the second and third order directional derivatives are:
\begin{subequations}
\begin{align}
\nabla^2 \varphi(w) [r] 
&\labrel{\ref{eq:Dinv}}{=} 
P(w^{-1}) r,
\label{eq:psid2}
\\
\nabla^3 \varphi(w) [r, r] 
&\labrel{\ref{eq:Pinv}}{=} 
-2 P(w^{-1/2}) (P(w^{-1/2}) r)^2.
\label{eq:psid3}
\end{align}
\end{subequations}
Note that unlike the separable spectral function case in \cref{sec:spectral:sep}, here we do not need the explicit eigenvalues of $w$.
\section{Cones and barrier functions}
\label{sec:cones}

In this paper we are concerned with a class of proper cones that can be characterized as follows:
\begin{equation}
\K \coloneqq \cl \{ \tu \in \cE : \zeta(\tu) \geq 0 \} \subset \tV,
\label{eq:K}
\end{equation}
where $\zeta : \cE \to \bbR$ is a concave, (degree one) homogeneous function and $\cE$ is some convex cone in the space $\tV$.
In particular, we define $\zeta$ in terms of a $C^3$-smooth spectral function $\varphi$ that is defined on the interior of a cone of squares $\cQ$ of a Jordan algebra $V$ of rank $d$.

\subsection{The homogeneous case}
\label{sec:cones:hom}

First, we suppose that $\varphi$ is convex and homogeneous.
Then $\zeta(u, w) \coloneqq u - \varphi(w)$ is concave and homogeneous, and we let $\cE \coloneqq \bbR \times \intr(\cQ)$ and $\tV \coloneqq \bbR \times V$.
This defines a convex cone that is the closure of the epigraph set of $\varphi$:
\begin{equation}
\K_h \coloneqq \cl \{ 
(u, w) \in \bbR \times \intr(\cQ) : 
u \geq \varphi (w) 
\}.
\label{eq:epi}
\end{equation}
Note that if $\varphi$ is concave, we can analogously define a cone from the hypograph set of $\varphi$.
In \cref{sec:rtdet}, we consider the root-determinant cone, which is the hypograph of the concave root-determinant function.
To check membership in $\intr(\K_h)$, we first determine whether $w \in \intr(\cQ)$ (which is equivalent to positivity of the eigenvalues), and if so, whether $\zeta(\tu) > 0$.

\subsection{The non-homogeneous case}
\label{sec:cones:per}

Now we suppose that $\varphi$ is convex and non-homogeneous.
We define the perspective function of $\varphi$, $\per \varphi : \bbR_{>} \times \intr(\cQ) \to \bbR$, as $(\per \varphi)(v, w) \coloneqq v \varphi(v^{-1} w)$.
This is a homogeneous and convex function \citep[Section 3.2.6]{boyd2004convex}.
We let $\zeta(u, v, w) \coloneqq u - (\per \varphi)(v, w)$, with $\cE \coloneqq \bbR \times \bbR_{>} \times \intr(\cQ)$ and $\tV \coloneqq \bbR \times \bbR \times V$.
This defines a convex cone that is the closure of the epigraph set of the perspective function of $\varphi$:
\begin{equation}
\K_p \coloneqq \cl \{ 
(u, v, w) \in \bbR \times \bbR_{>} \times \intr(\cQ) : 
u \geq v \varphi (v^{-1} w) 
\}.
\label{eq:epiper}
\end{equation}
Equivalently, we can view $\K_p$ as the closed conic hull of the epigraph set of $\varphi$ \citep[Chapter 5]{nesterov1994interior}.
In \cref{sec:matmono}, we consider the special case where $\varphi$ is a separable spectral function with matrix monotone first derivative.
The membership check for $\intr(\K_p)$ is similar to that of $\intr(\K_h)$ except we first check whether $v > 0$.

\subsection{Dual cones}
\label{sec:cones:dual}

A proper cone is \emph{primitive} if it cannot be written as a Cartesian product of two or more lower dimensional proper cones.
The dual cone of a primitive, proper cone $\K$ is another primitive, proper cone:
\begin{equation}
\K^\ast \coloneqq \{ z : \langle s, z \rangle \geq 0, \forall s \in \K \}.
\end{equation}
Recall that when $\K$ is defined through Hypatia's generic cone interface, both $\K$ and $\K^\ast$ become available for constructing conic models. 

We assume that $\varphi$ is convex, and we derive the dual cones of the epigraph cones $\K_h$ and $\K_p$ (these steps can be adapted for analogous hypograph cones if $\varphi$ is concave).
We define the convex conjugate function $\varphi^\ast : V \to \bbR \cup \infty$ of $\varphi$ as the modified Legendre-Fenchel transformation (similar to \citet[page 483]{zhang2004new}):
\begin{equation}
\varphi^\ast (r) = \textstyle\sup_{w \in \dom(\varphi)} 
\{ -\langle w, r \rangle - \varphi(w) \},
\label{eq:conj}
\end{equation}
which is a convex function.
The conjugate of a symmetric function is also a symmetric function \citep[Lemma 29]{baes2007convexity}, and the conjugate of a spectral function induced by a symmetric function $f$ is the spectral function induced by $f^\ast$ \citep[Theorem 30]{baes2007convexity}.
Thus for $\varphi(w) = f(\lambda(w))$ we have the conjugate function $\varphi^\ast(w) = f^\ast(\lambda(w))$.

For the epigraph-perspective cone $\K_p$ in \cref{eq:epiper}, \citet[Theorem 3.2]{zhang2004new} and \citet[Theorem 14.4]{rockafellar2015convex} derive the dual cone $\K_p^\ast$:
\begin{equation}
\K_p^\ast = \cl \{ 
(u, v, w) \in \bbR_{>} \times \bbR \times V : 
v \geq u \varphi^\ast (u^{-1} w) 
\}.
\label{eq:epiperdual}
\end{equation}
We can view $\K_p^\ast$ as the epigraph set of the perspective function of the conjugate of $\varphi$, but with the epigraph and perspective components swapped (compare to \cref{eq:epiper}).
Depending on the natural domain of $\varphi^\ast$, the $w$ component of $\K_p^\ast$ is not necessarily restricted to lie in $\cQ$; in \cref{sec:matmono:cone} we discuss several example spectral functions, some of which have conjugates defined on all $V$ and others only on $\intr(\cQ)$.

For $\K_h$ in \cref{eq:epi}, we derive the dual cone $\K_h^\ast$ as follows.
Since $\varphi$ is homogeneous in this case, $(\per \varphi)(v, w) = v \varphi (v^{-1} w) = \varphi(w)$.
Therefore the corresponding perspective cone $\K_p$ for $\varphi$ is not a primitive cone, as it can be written as a (permuted) Cartesian product of $\Rge$ and $\K_h$:
\begin{equation}
\K_p = \cl \{ 
(u, v, w) \in \tV : v \in \Rge, (u, w) \in \K_h
\}.
\end{equation}
Since the dual cone of a Cartesian product of cones is the Cartesian product of their dual cones, we have (since $\Rge^\ast = \Rge$):
\begin{equation}
\K_p^\ast = \cl \{ 
(u, v, w) \in \tV : v \in \Rge, (u, w) \in \K_h^\ast
\}.
\label{eq:homK}
\end{equation}
By \citep[Theorem 2.1]{lasserre1998homogeneous}, the homogeneity of $\varphi$ implies that $\varphi^\ast$ can only take the values zero or infinity.
Hence by \cref{eq:epiperdual}, we know:
\begin{equation}
\K_p^\ast = \cl \{ 
(u, v, w) \in \bbR_{>} \times \bbR \times V : 
v \geq 0, u \varphi^\ast (u^{-1} w) < \infty
\}.
\label{eq:homK2}
\end{equation}
Since \cref{eq:homK,eq:homK2} describe the same cone, we can conclude that the dual cone of $\K_h$ is:
\begin{equation}
\K_h^\ast = \cl \{ 
(u, w) \in \bbR_{>} \times V : 
\varphi^\ast (u^{-1} w) < \infty \}.
\label{eq:epidualhomog}
\end{equation}

\subsection{Barrier functions and oracles}
\label{sec:cones:bar}

A logarithmically homogeneous barrier (LHB) function $\Gamma$ for a proper cone $\K \subset \tV$ is $C^2$-smooth and satisfies $\Gamma(\tu_i) \to \infty$ along every sequence $\tu_i \in \intr(\K)$ converging to the boundary of $\K$, and:
\begin{align}
\Gamma(\theta \tu) 
& = \Gamma(\tu) - \nu \log(\theta) 
& \forall \tu \in \K, \theta \in \bbR,
\label{eq:lhscb:lh}
\end{align}
for some $\nu \geq 0$ \citep[Definition 2.3.2]{nesterov1994interior}.
If $\Gamma$ is also self-concordant, then it is an LHSCB with parameter $\nu \geq 1$ (or a $\nu$-LHSCB) for $\K$.
For self-concordance, $\Gamma$ must be $C^3$-smooth and satisfy \citep[Definition 2.1.1]{nesterov1994interior}:
\begin{align}
\lvert \nabla^3 \Gamma(\tu) [\tp, \tp, \tp] \rvert 
& \leq 2 ( \nabla^2 \Gamma(\tu) [\tp, \tp] )^{3/2} 
& \forall \tu \in \intr(\K), \tp \in \tV.
\label{eq:lhscb:sc}
\end{align}
The best known interior point algorithms need at most $O(\sqrt{\nu} \log (1 / \varepsilon))$ iterations to converge to a solution within $\varepsilon$ tolerance \citep{nesterov1997self}.

The LHB we consider for a cone of the form \cref{eq:K} is:
\begin{equation}
\Gamma(\tu) \coloneqq -\log(\zeta(\tu)) + \Psi(\tu),
\label{eq:lhb}
\end{equation}
where $\Psi$ can be thought of as an LHSCB for the domain of $\zeta$ or $\cl(\cE)$. 
The negative logarithm function $-\log$ is the standard LHSCB for $\Rge$, with parameter $\nu = 1$.
Similarly, the spectral function $-\logdet$ (see \cref{sec:spectral:log}) is the standard LHSCB for a cone of squares $\cQ$ of $V$, with $\nu = d$ (the rank of $V$).
For $\K_h$ we let $\Psi(\tu) = -\logdet(w)$, hence $\Gamma$ has parameter $\nu = 1 + d$.
Since an LHSCB for a Cartesian product of cones is the sum of LHSCBs for the primitive cones, for $\K_p$ we let $\Psi(\tu) = -\log(v) - \logdet(w)$, hence $\Gamma$ has parameter $\nu = 2 + d$.
Note that although $\Gamma$ is an LHB, it is not necessarily self-concordant; in \cref{sec:matmono:lhscb,sec:rtdet:bar} we prove that $\Gamma$ is an LHSCB for some useful special cases.

We now define four barrier oracles that Hypatia's PDIPM uses; for ideal performance, these oracle implementations should be efficient and numerically stable.
For an interior point $\tu \in \intr(\K)$ and a direction $\tp \in \tV$, the gradient $g$, the Hessian product $H$, the inverse Hessian product $\bH$, and the third order directional derivative $T$ are:
\begin{subequations}
\begin{align}
g &\coloneqq \nabla \Gamma(\tu),
\label{eq:oracles:g}
\\
H &\coloneqq \nabla^2 \Gamma(\tu) [\tp],
\\
\bH &\coloneqq (\nabla^2 \Gamma(\tu))^{-1} [\tp].
\label{eq:oracles:bh}
\\
T &\coloneqq \nabla^3 \Gamma(\tu) [\tp, \tp].
\end{align}
\label{eq:oracles}
\end{subequations}
Note $g, H, \bH, T \in \tV$. 
In later sections, we use subscripts to refer to subcomponents of these oracles, for example the $w$ component of the gradient oracle is $g_w \coloneqq \nabla_w \Gamma(\tu) \in V$.
Ideally, $H$ applies the positive definite linear operator $\nabla^2 \Gamma (\tu) : \tV \to \tV$ without constructing an explicit Hessian, and similarly, $\bH$ applies the (unique) inverse operator $(\nabla^2 \Gamma (\tu))^{-1} : \tV \to \tV$ without constructing or factorizing an explicit Hessian.

We note that for the standard LHSCB $\Psi$ for a cone of squares, efficient and numerically stable procedures for these four oracles are well-known.
The same cannot be said for the LHB $\Gamma$ currently.
In \cref{sec:epiper}, we derive these oracles for $\K_p$ (noting that they can be adapted easily for $\K_h$).
In the special cases for which we show $\Gamma$ is an LHSCB, the oracles can be computed particularly efficiently.
\section{Barrier oracles for epigraph-perspective cones}
\label{sec:epiper}

We consider the epigraph-perspective cone $\K_p$ defined in \cref{eq:epiper}.
Recall that we let $\tp = (p, q, r) \in \bbR \times \bbR \times V$ and $\tu = (u, v, w) \in \intr(\K_p)$, and we define $\zeta$ and $\Gamma : \intr(\K_p) \to \bbR$ from \cref{eq:lhb} as:
\begin{subequations}
\begin{align}
\zeta (\tu) &\coloneqq u - v \varphi(v^{-1} w),
\\
\Gamma (\tu) &\coloneqq -\log(\zeta(\tu)) - \log(v) - \logdet(w).
\end{align}
\end{subequations}
In this section, we derive expressions and evaluation procedures for the $g$, $H$, $T$, and $\bH$ oracles (defined in \cref{sec:cones:bar}) corresponding to the LHB $\Gamma$ for $\K_p$.
We note that the oracles for $\K_h$ in \cref{eq:epi} are simpler because no perspective operation is needed for a homogeneous $\varphi$; they can be obtained by fixing $v = 1$ and $q = 0$ and ignoring the $v$ components in the oracle expressions in this section.

Without assuming any particular form for $\varphi$, we write $g$, $H$, and $T$ in \cref{sec:epiper:deriv} and $\bH$ in \cref{sec:epiper:invhessgen} in terms of the derivatives of $\varphi$.
If $\varphi$ is a spectral function, these derivatives can be computed using the expressions from \cref{sec:spectral}.
In the case that $\varphi$ is a separable spectral function (see \cref{sec:spectral:sep}), we derive a more specialized procedure for $\bH$ in \cref{sec:epiper:invhesssep}, which is no more expensive than $H$.
Finally, in \cref{sec:epiper:logdet}, we specialize the four oracles for the negative log-determinant function (i.e. $\varphi(w) = -\logdet(w)$; see \cref{sec:spectral:log}) and we discuss implementations.

\subsection{Derivatives}
\label{sec:epiper:deriv}

First, we express the derivatives of $\zeta$ in terms of those of $\varphi$.
We define the function $\mu : \bbR_{>} \times \intr(\cQ) \to \intr(\cQ)$ and its first directional derivative $\xi \in V$ in the direction $(q, r)$ as:
\begin{subequations}
\begin{align}
\mu(v, w) &\coloneqq v^{-1} w,
\\
\xi &\coloneqq \nabla \mu(v, w) [(q, r)]
= \nabla_v \mu(v, w) q + \nabla_w \mu(v, w) [r]
= v^{-1} (r - q \mu(v, w)).
\end{align}
\label{eq:mu}
\end{subequations}
For convenience, we fix the constants $\mu \coloneqq \mu(v, w)$, $\varphi \coloneqq \varphi(\mu)$, and $\zeta \coloneqq \zeta(\tu)$.
Let $\nabla \varphi$, $\nabla^2 \varphi$, and $\nabla^3 \varphi$ be the derivatives of $\varphi$ evaluated at $\mu$, and let $\nabla \zeta$, $\nabla^2 \zeta$, and $\nabla^3 \zeta$ be the derivatives of $\zeta$ evaluated at $\tu$.
Using \cref{eq:mu}, the directional derivatives of $\zeta$ can be written compactly as:
\begin{subequations}
\begin{align}
\nabla_u \zeta
&= 1,
\\
\nabla_v \zeta
&= -\varphi + \nabla \varphi [\mu],
\\
\nabla_w \zeta
&= -\nabla \varphi,
\\
\nabla \zeta [\tp] 
&= p - q \varphi - v \nabla \varphi [\xi],
\\
(\nabla^2 \zeta [\tp])_v
&= \nabla^2 \varphi [\xi, \mu],
\\
(\nabla^2 \zeta [\tp])_w
&= -\nabla^2 \varphi [\xi],
\\
\nabla^2 \zeta [\tp, \tp] 
&= -v \nabla^2 \varphi [\xi, \xi],
\\
(\nabla^3 \zeta [\tp, \tp])_v
&= \nabla^3 \varphi [\xi, \xi, \mu] + \nabla^2 \varphi [\xi, \xi] - 
2 v^{-1} q \nabla^2 \varphi [\xi, \mu],
\\
(\nabla^3 \zeta [\tp, \tp])_w
&= 2 v^{-1} q \nabla^2 \varphi [\xi] - \nabla^3 \varphi [\xi, \xi],
\\
\nabla^3 \zeta [\tp, \tp, \tp] 
&= -v \nabla^3 \varphi [\xi, \xi, \xi] + 3 q \nabla^2 \varphi [\xi, \xi].
\end{align}
\label{eq:zetadds}
\end{subequations}

Using \cref{eq:zetadds}, we now derive the directional derivatives of $\Gamma$.
For convenience, we let $\nabla \Gamma$, $\nabla^2 \Gamma$, $\nabla^3 \Gamma$ be the derivatives of $\Gamma$ evaluated at $\tu$.
We define:
\begin{equation}
\sigma \coloneqq -\nabla_v \zeta = \varphi - \nabla \varphi[\mu]
\in \bbR.
\label{eq:epiper:sigma}
\end{equation}
The components of the gradient $g$ of $\Gamma$ are:
\begin{subequations}
\begin{align}
g_u
&= -\zeta^{-1},
\\
g_v
&= \zeta^{-1} \sigma - v^{-1},
\\
g_w
&\labrel{\ref{eq:Dpsi}}{=} \zeta^{-1} \nabla \varphi - w^{-1}.
\label{eq:Dwgamma}
\end{align}
\label{eq:epipergrad}
\end{subequations}
Differentiating \cref{eq:epipergrad}, the Hessian components are:
\begin{subequations}
\begin{align}
\nabla^2_{u,u} \Gamma
&= \zeta^{-2} 
> 0,
\\
\nabla^2_{v,u} \Gamma
&= -\zeta^{-2} \sigma 
\in \bbR,
\\
\nabla^2_{w,u} \Gamma
&= -\zeta^{-2} \nabla \varphi 
\in V,
\\
\nabla^2_{v,v} \Gamma
&= v^{-2} + \zeta^{-2} \sigma^2 + v^{-1} \zeta^{-1} \nabla^2 \varphi [\mu, \mu]
> 0,
\\
\nabla^2_{w,v} \Gamma
&= \zeta^{-2} \sigma \nabla \varphi - v^{-1} \zeta^{-1} \nabla^2 \varphi [\mu]
\in V.
\end{align}
\label{eq:D2uvw}
\end{subequations}
Differentiating \cref{eq:Dwgamma} in the direction $r$:
\begin{align}
\nabla^2_{w,w} \Gamma [r] 
&\labrel{\ref{eq:Dinv}}{=} \zeta^{-2} \nabla \varphi [r] \nabla \varphi + 
v^{-1} \zeta^{-1} \nabla^2 \varphi [r] + P(w^{-1}) r
\in V.
\label{eq:D2ww}
\end{align}
Let:
\begin{equation}
\chi \coloneqq \zeta^{-1} (p - q \sigma - \nabla \varphi[r]) \in \bbR.
\label{eq:epiper:chi}
\end{equation}
The components of the Hessian product $H$ are:
\begin{subequations}
\begin{align}
H_u 
&= \zeta^{-1} \chi,
\\
H_v 
&= -\zeta^{-1} \sigma \chi - \zeta^{-1} \nabla^2 \varphi [\xi, \mu] + v^{-2} q,
\\
H_w 
&= -\zeta^{-1} \chi \nabla \varphi + \zeta^{-1} \nabla^2 \varphi [\xi] + P(w^{-1}) r.
\end{align}
\label{eq:hessprod}
\end{subequations}
Let:
\begin{equation}
\kappa \coloneqq
2 \zeta^{-1} (\chi + v^{-1} q) \nabla^2 \varphi [\xi] - 
\zeta^{-1} \nabla^3 \varphi [\xi, \xi]
\in V.
\label{eq:epiper:kappa}
\end{equation}
The components of the third order directional derivative $T$ are:
\begin{subequations}
\begin{align}
T_u 
&= -2 \zeta^{-1} \chi^2 - v \zeta^{-2} \nabla^2 \varphi [\xi, \xi],
\\
T_v 
&= -T_u \sigma + \langle \kappa, \mu \rangle - \zeta^{-1} \nabla^2 \varphi [\xi, \xi] - 2 q^2 v^{-3},
\\
T_w 
&\labrel{\ref{eq:Pinv}}{=} 
-T_u \nabla \varphi - \kappa - 2 P(w^{-1/2}) (P(w^{-1/2}) r)^2.
\end{align}
\label{eq:third}
\end{subequations}

\subsection{Inverse Hessian operator}
\label{sec:epiper:invhessgen}

The Hessian of $\Gamma$ at any point $\tu \in \intr(\K_p)$ is a positive definite linear operator and hence invertible.
By treating the components of the Hessian in \cref{eq:D2uvw,eq:D2ww} analogously to blocks of a positive definite matrix, we derive the inverse operator.
For convenience, we let:
\begin{subequations}
\begin{align}
Y_u &\coloneqq (\nabla^2_{w,w} \Gamma)^{-1} \nabla^2_{w,u} \Gamma,
\\
Y_v &\coloneqq (\nabla^2_{w,w} \Gamma)^{-1} \nabla^2_{w,v} \Gamma,
\\
Z_{u,u} &\coloneqq \nabla^2_{u,u} \Gamma - \langle \nabla^2_{w,u} \Gamma, Y_u \rangle,
\\
Z_{v,u} &\coloneqq \nabla^2_{v,u} \Gamma - \langle \nabla^2_{w,u} \Gamma, Y_v \rangle,
\\
Z_{v,v} &\coloneqq \nabla^2_{v,v} \Gamma - \langle \nabla^2_{w,v} \Gamma, Y_v \rangle.
\end{align}
\label{eq:defYZ}
\end{subequations}
Note $Y_u, Y_v \in V$.
We let $Z$ be:
\begin{equation}
Z \coloneqq 
\begin{bmatrix}
Z_{u,u} & Z_{v,u} \\
Z_{v,u} & Z_{v,v}
\end{bmatrix}
\in \bbS_{\succ}^2,
\label{eq:Z}
\end{equation}
and its inverse is:
\begin{equation}
\bZ \coloneqq Z^{-1}
= \frac{1}{Z_{u,u} Z_{v,v} - Z_{v,u}^2}
\begin{bmatrix}
Z_{v,v} & -Z_{v,u} \\
-Z_{v,u} & Z_{u,u}
\end{bmatrix}
\in \bbS_{\succ}^2.
\label{eq:barZ}
\end{equation}
It can be verified (for example, by analogy to the block symmetric matrix inverse formula) that the inverse Hessian product oracle $\bH$ in \cref{eq:oracles:bh} is:
\begin{subequations}
\begin{align}
\bH_u 
&= \bZ_{u,u} (p - \langle Y_u, r \rangle ) + \bZ_{v,u} (q - \langle Y_v, r \rangle ),
\\
\bH_v 
&= \bZ_{v,u} (p - \langle Y_u, r \rangle ) + \bZ_{v,v} (q - \langle Y_v, r \rangle ),
\\
\bH_w 
&= -\bH_u Y_u - \bH_v Y_v + (\nabla^2_{w,w} \Gamma)^{-1} r
.
\end{align}
\label{eq:general:invprod}
\end{subequations}
Hence computing $\bH$ is essentially only as difficult as applying the positive definite linear operator $(\nabla^2_{w,w} \Gamma)^{-1}$.
We are not aware of a simple expression for $(\nabla^2_{w,w} \Gamma)^{-1}$ in general, but we explore the special cases of separable spectral functions in \cref{sec:epiper:invhesssep}, the negative log-determinant function in \cref{sec:epiper:logdet}, and the root-determinant function in \cref{sec:rtdet:oracles}.

\subsection{Inverse Hessian operator for the separable spectral case}
\label{sec:epiper:invhesssep}


Suppose $w \succ 0$ has the spectral decomposition \cref{eq:specdecomp}, i.e. $w$ has the eigenvalues $w_1, \ldots, w_d > 0$ and the Jordan frame $c_1, \ldots, c_d$.
As in \cref{sec:spectral:sep}, we assume distinct eigenvalues for simplicity.
In the special case where $\varphi$ is a convex separable spectral function, i.e. $\varphi(w) = \sum_{i \in \iin{d}} h(w_i)$ for some convex $h : \bbR_{>} \to \bbR$, we show how to compute $\bH$ as efficiently as Hessian product oracle $H$.
For all $i \in \iin{d}$, we let $h_i$, $(\nabla h)_i$, $(\nabla^2 h)_i$, and $(\nabla^3 h)_i$ denote the value and derivatives of $h$ evaluated at $\mu$.
We define $m_{i,j} \in \bbR$ for $i, j \in \iin{d}$ as:
\begin{align}
m_{i,j} &\coloneqq 
\begin{cases}
\zeta^{-1} \frac{ (\nabla h)_i - (\nabla h)_j}{w_i - w_j} + 
w_i^{-1} w_j^{-1} & i \neq j,
\\
\zeta^{-1} v^{-1} (\nabla^2 h)_i + w_i^{-2} & i = j
.
\end{cases}
\end{align}
Since $h$ is convex, $m_{i,j} > 0, \forall i, j \in \iin{d}$.
Let $M: V \to V$ be the self-adjoint linear operator:
\begin{align}
M &\coloneqq 
v^{-1} \zeta^{-1} \nabla^2 \varphi + P(w^{-1})
= \sum_{i,j \in \iin{d} : i < j} 4 m_{i,j} L(c_i) L(c_j) + \sum_{i \in \iin{d}} m_{i,i} P(c_i).
\label{eq:defM}
\end{align}
Using \cref{eq:opinv}, we have the self-adjoint inverse operator of $M$:
\begin{align}
M^{-1} 
&= \sum_{i,j \in \iin{d} : i < j} 4 m_{i,j}^{-1} L(c_i) L(c_j) + 
\sum_{i \in \iin{d}} m_{i,i}^{-1} P(c_i).
\label{eq:defMi}
\end{align}
Substituting \cref{eq:defM} into \cref{eq:D2ww}, we have:
\begin{align}
\nabla^2_{w,w} \Gamma [r] 
&= \zeta^{-2} \nabla \varphi [r] \nabla \varphi + M r.
\label{eq:D2wwalg}
\end{align}
Note that the first term in \cref{eq:D2wwalg} is analogous to the application (to $r$) of a low-rank update to $M$, and that $M^{-1}$ in \cref{eq:defMi} is easy to apply.
By analogy to the Sherman-Morrison-Woodbury formula \citep[Theorem 1.1]{deng2011generalization}, we can derive a simple expression for the inverse operator $(\nabla_{w,w}^2 \Gamma)^{-1} r$.

We let:
\begin{subequations}
\begin{align}
\alpha 
&\coloneqq M^{-1} \nabla \varphi
= \sum_{i \in \iin{d}} m_{i,i}^{-1} (\nabla h)_i c_i
\in V,
\\
\gamma 
&\coloneqq v^{-2} \zeta^{-1} M^{-1} \nabla^2 \varphi[w]
= v^{-2} \zeta^{-1} \sum_{i \in \iin{d}} m_{i,i}^{-1} (\nabla^2 h)_i w_i c_i
\succ 0.
\end{align}
\label{eq:greekconsts}
\end{subequations}
Noting that $\gamma, w^{-1} \succ 0$ implies $\langle \gamma, w^{-1} \rangle > 0$, we define the scalar constants:
\begin{subequations}
\begin{align}
k_1 
&\coloneqq \zeta^2 + \langle \nabla \varphi, \alpha \rangle 
> 0,
\\
k_2 
&\coloneqq \sigma + \langle \nabla \varphi, \gamma \rangle
= \sigma + v^{-2} \zeta^{-1} \langle \nabla^2 \varphi[w], \alpha \rangle,
\\
k_3 &\coloneqq 
v^{-2} + v^{-2} \zeta^{-1} \langle \nabla^2 \varphi[w], \mu - \gamma \rangle
= v^{-2} + v^{-1} \langle \gamma, w^{-1} \rangle
> 0.
\label{eq:defk3}
\end{align}
\label{eq:kconsts}
\end{subequations}
Now using the Sherman-Morrison-Woodbury formula:
\begin{subequations}
\begin{align}
(\nabla^2_{w,w} \Gamma)^{-1} r 
&= M^{-1} r - \frac{
\zeta^{-2} \langle M^{-1} \nabla \varphi, r \rangle
}{
1 + \zeta^{-2} \langle M^{-1} \nabla \varphi, \nabla \varphi \rangle
} M^{-1} \nabla \varphi
\\
&= M^{-1} r - k_1^{-1} \langle \alpha, r \rangle \alpha.
\end{align}
\label{eq:sm}
\end{subequations}
Substituting \cref{eq:D2uvw,eq:sm} into \cref{eq:defYZ}:
\begin{subequations}
\begin{align}
Y_u 
&= (\nabla^2_{w,w} \Gamma)^{-1} (-\zeta^{-2} \nabla \varphi) 
\\
&= -\zeta^{-2} \alpha + \zeta^{-2} k_1^{-1} \langle \alpha, \nabla \varphi \rangle \alpha
\\
&= -k_1^{-1} \alpha,
\label{eq:Yu}
\\
Y_v 
&= (\nabla^2_{w,w} \Gamma)^{-1} (\zeta^{-2} \sigma \nabla \varphi - 
v^{-1} \zeta^{-1} \nabla^2 \varphi [\mu]) 
\\
&= -\sigma Y_u - v^{-2} \zeta^{-1} (\nabla^2_{w,w} \Gamma)^{-1} \nabla^2 \varphi [w]
\\
&= \sigma k_1^{-1} \alpha - \gamma + v^{-2} \zeta^{-1} k_1^{-1} \langle \alpha, \nabla^2 \varphi[w] \rangle \alpha
\\
&= k_1^{-1} k_2 \alpha - \gamma,
\\
Z_{u,u} 
&= \zeta^{-2} - \langle \nabla^2_{w,u} \Gamma, Y_u \rangle
\label{eq:Zuu}
\\
&= \zeta^{-2} - \zeta^{-2} k_1^{-1} \langle \nabla \varphi, \alpha \rangle
\\
&= k_1^{-1},
\\
Z_{v,u} 
&= -\zeta^{-2} \sigma - \langle \nabla^2_{w,u} \Gamma, Y_v \rangle
\\
&= -\zeta^{-2} (\sigma - \langle \nabla \varphi, k_1^{-1} k_2 \alpha - \gamma \rangle)
\\
&= -\zeta^{-2} (\sigma - k_1^{-1} k_2 (k_1 - \zeta^2) + k_2 - \sigma)
\\
&= -k_1^{-1} k_2,
\label{eq:Zvu}
\\ 
Z_{v,v} 
&= \nabla^2_{v,v} \Gamma - \langle \nabla^2_{w,v} \Gamma, Y_v \rangle
\\
&= \nabla^2_{v,v} \Gamma + 
\sigma \langle \nabla^2_{w,u} \Gamma, Y_v \rangle + 
v^{-2} \zeta^{-1} \langle \nabla^2 \varphi [w], Y_v \rangle
\\
&\labrel{\ref{eq:Zvu}}{=} 
\nabla^2_{v,v} \Gamma +
\sigma (k_1^{-1} k_2 - \zeta^{-2} \sigma) + 
v^{-2} \zeta^{-1} \langle \nabla^2 \varphi [w], Y_v \rangle
\\
&= v^{-2} + v^{-3} \zeta^{-1} \nabla^2 \varphi [w, w] + 
\sigma k_1^{-1} k_2 + 
v^{-2} \zeta^{-1} \langle \nabla^2 \varphi [w], Y_v \rangle
\\
&= v^{-2} + v^{-3} \zeta^{-1} \nabla^2 \varphi [w, w] + 
\sigma k_1^{-1} k_2 + 
k_1^{-1} k_2 (k_2 - \sigma) - 
v^{-2} \zeta^{-1} \langle \nabla^2 \varphi [w], \gamma \rangle
\\
&= v^{-2} + k_1^{-1} k_2^2 +
v^{-2} \zeta^{-1} \langle \nabla^2 \varphi [w], \mu - \gamma \rangle
\\
&= k_3 + k_1^{-1} k_2^2.
\end{align}
\label{eq:YZ}
\end{subequations}
For $Z$ in \cref{eq:Z}, we have $\det(Z) = k_1^{-1} k_3$, so its inverse $\bZ$ in \cref{eq:barZ} is:
\begin{subequations}
\begin{align}
\bZ_{u,u} &= k_1 (k_3 + k_1^{-1} k_2^2) k_3^{-1} 
= k_1 + k_2^2 k_3^{-1},
\\
\bZ_{u,v} &= k_2 k_3^{-1},
\\
\bZ_{v,v} &= k_3^{-1}.
\end{align}
\label{eq:barZsep}
\end{subequations}
Finally, we substitute \cref{eq:YZ,eq:barZsep} into \cref{eq:general:invprod} to derive the inverse Hessian product $\bH$.
We let:
\begin{subequations}
\begin{align}
c_1 &\coloneqq p - \langle Y_u, r \rangle 
\labrel{\ref{eq:YZ}}{=} 
p + k_1^{-1} \langle \alpha, r \rangle,
\\
c_2 &\coloneqq q - \langle Y_v, r \rangle 
\labrel{\ref{eq:YZ}}{=} 
q - k_1^{-1} k_2 \langle \alpha, r \rangle + \langle \gamma, r \rangle.
\end{align}
\end{subequations}
For convenience, we derive $\bH_v$ before $\bH_u$ and $\bH_w$:
\begin{subequations}
\begin{align}
\bH_v &= 
\bZ_{u,v} c_1 + \bZ_{v,v} c_2
\\
&= k_3^{-1} (k_2 c_1 + c_2)
\label{eq:bHvc2}
\\
&= k_3^{-1} (k_2 p + q + \langle \gamma, r \rangle),
\\
\bH_u
&= \bZ_{u,u} c_1 + \bZ_{u,v} c_2
\\
&\labrel{\ref{eq:bHvc2}}{=} 
(\bZ_{u,u} - \bZ_{u,v} k_2) c_1 + \bZ_{u,v} k_3 \bH_v
\\
&= k_1 p + k_2 \bH_v + \langle \alpha, r \rangle,
\\
\bH_w
&= -\bH_u Y_u - \bH_v Y_v + (\nabla^2_{w,w})^{-1} r
\\
&\labrel{\ref{eq:sm}}{=} 
\bH_u k_1^{-1} \alpha - \bH_v (k_1^{-1} k_2 \alpha - \gamma) + 
M^{-1} r - k_1^{-1} \langle \alpha, r \rangle \alpha
\\
&= p \alpha + \bH_v \gamma + M^{-1} r.
\end{align}
\label{eq:invhessprodsep}
\end{subequations}

In \cref{sec:testing:invhessprod}, we compare the efficiency and numerical performance of the closed-form formula for $\bH$ in \cref{eq:invhessprodsep} against a naive approach to computing $\bH$ that performs a Cholesky factorization of an explicit Hessian matrix and uses a direct linear solve.
The closed-form formula is faster and more scalable, more memory-efficient, more reliable to compute (as the Cholesky decomposition can fail), and more numerically accurate.

\subsection{Oracles for the log-determinant case}
\label{sec:epiper:logdet}

We now specialize the oracles derived in \cref{sec:epiper:deriv,sec:epiper:invhesssep} for the separable spectral function $\varphi(w) = -\logdet(w) = -\sum_{i \in \iin{d}} \log(w_i)$. 
In \cref{sec:matmono}, we show that $\Gamma$ is an LHSCB in this case.
We let:
\begin{equation}
\hxi \coloneqq P(w^{-1/2}) \xi
= v^{-1} P(w^{-1/2}) (-v^{-1} q w + r)
= v^{-1} (-v^{-1} q e + \hr) 
\in V.
\end{equation}
We have:
\begin{subequations}
\begin{align}
\nabla \varphi 
&\labrel{\ref{eq:Dpsi}}{=}
-\mu^{-1}
= -v P(w^{-1/2}) e,
\\
\nabla^2 \varphi [\xi] 
&\labrel{\ref{eq:psid2}}{=} 
v^2 P(w^{-1}) \xi
= v^2 P(w^{-1/2}) \hxi,
\label{eq:log:nabla2}
\\
\nabla^3 \varphi [\xi, \xi] 
&\labrel{\ref{eq:psid3}}{=} 
-2 v^3 P(w^{-1/2}) \hxi^2.
\label{eq:log:nabla3}
\end{align}
\label{eq:logDs}
\end{subequations}
The constants from \cref{sec:epiper:deriv} have the form:
\begin{subequations}
\begin{align}
\sigma 
&\labrel{\ref{eq:epiper:sigma}}{=}
\varphi + d,
\\
\chi 
&\labrel{\ref{eq:epiper:chi}}{=}
\zeta^{-1} (p - q \sigma + v \tr(\hr)),
\end{align}
\end{subequations}

From \cref{eq:epipergrad}, the $w$ component of the gradient is:
\begin{subequations}
\begin{align}
g_w
&= -(1 + v \zeta^{-1}) w^{-1}.
\end{align}
\end{subequations}
From \cref{eq:hessprod}, the $v$ and $w$ components of the Hessian product are:
\begin{subequations}
\begin{align}
H_v
&= -\zeta^{-1} \sigma \chi - v \zeta^{-1} \tr(\hxi) + v^{-2} q,
\\
H_w
&= P(w^{-1/2}) ( v \zeta^{-1} \chi e + v^2 \zeta^{-1} \hxi + \hr ).
\end{align}
\end{subequations}
From \cref{eq:third}, the third order directional derivative is:
\begin{subequations}
\begin{align}
T_u
&= -2 \zeta^{-1} \chi^2 - v^3 \zeta^{-2} \tr(\hxi^2),
\\
T_v
&= -T_u \sigma +
v \zeta^{-1} ( 2 (\chi + v^{-1} q) \tr(\hxi) + v \tr(\hxi^2) ) - 2 v^{-3} q^2,
\\
T_w
&= P(w^{-1/2}) (T_u v e - 2 \zeta^{-1} v^2 ((\chi + v^{-1} q) \hxi + v \hxi^2) - 2 \hr^2).
\end{align}
\end{subequations}


We derive the inverse Hessian product $\bH$ by specializing the separable case in \cref{eq:invhessprodsep}.
We let:
\begin{subequations}
\begin{align}
\chr &\coloneqq P(w^{1/2}) r \in V,
\\
\theta &\coloneqq v^2 (\zeta + (1 + d) v)^{-1}.
\end{align}
\end{subequations}
From \cref{eq:greekconsts,eq:kconsts}, we have:
\begin{subequations}
\begin{align}
M^{-1} &= \zeta (\zeta + v)^{-1} P(w),
\\
\alpha &= -v \zeta (\zeta + v)^{-1} w,
\\
\gamma &= (\zeta + v)^{-1} w,
\\
k_1 &= \zeta^2 + d v^2 \zeta (\zeta + v)^{-1},
\\
k_2 &= \varphi + d \zeta (\zeta + v)^{-1},
\\
k_3^{-1} &= (\zeta + v) \theta.
\end{align}
\end{subequations} 
For convenience, we derive $\bH_v$ before $\bH_u$ and $\bH_w$ as in \cref{eq:invhessprodsep}:
\begin{subequations}
\begin{align}
\bH_v 
&= k_3^{-1} (k_2 p + q + \langle \gamma, r \rangle)
\\
&= (\zeta + v) \theta ((\varphi + d \zeta (\zeta + v)^{-1}) p + 
q + \langle r, w \rangle (\zeta + v)^{-1})
\\
&= \theta ( (\zeta + v) (\varphi p + q) + d \zeta p + \tr(\chr) ).
\\
\bH_u 
&= k_1 p + (\varphi + d \zeta (\zeta + v)^{-1}) \bH_v + \langle \alpha, r \rangle
\\
&= (\zeta^2 + d v^2 \zeta (\zeta + v)^{-1}) p + (\varphi + d \zeta (\zeta + v)^{-1}) \bH_v + 
\langle -v \zeta (\zeta + v)^{-1} w, r \rangle
\\
&= \zeta (\zeta + v)^{-1} (d v^2 p + d \bH_v - v \langle w, r \rangle) + \zeta^2 p + \varphi \bH_v.
\\
\bH_w 
&= p \alpha + \bH_v \gamma + M^{-1} r
\\
&= -p v \zeta (\zeta + v)^{-1} w + \bH_v (\zeta + v)^{-1} w + \zeta (\zeta + v)^{-1} P(w) r
\\
&= (\zeta + v)^{-1} P(w^{1/2}) ((-\zeta p v + \bH_v)e + \zeta \chr).
\end{align}
\label{eq:log:invprod}
\end{subequations}

Recall that in \cref{sec:epiper:invhesssep}, we used the simplifying assumption of distinct eigenvalues, but for the negative log-determinant case this is not necessary.
Note that if it is possible to apply $P(w^{1/2})$ and $P(w^{-1/2})$ without accessing the eigenvalues of $w$, then all four oracles can be computed without an explicit eigendecomposition.
For example, in our implementation for $V = \bbS^d$ and $V = \bbH^d$, only a Cholesky factorization of $w$ is needed.
This is unlike the more general separable spectral function case, where the explicit eigenvalues of $w$ are needed.

\section{Matrix monotone derivative cones}
\label{sec:matmono}

After defining the matrix monotone property of a function in \cref{sec:matmono:mm}, we introduce the matrix monotone derivative cone $\kmmd$ in \cref{sec:matmono:cone}.
$\kmmd$ is a special case of the epigraph-perspective cone $\K_p$ with a separable spectral function $\varphi$.
In \cref{sec:matmono:lhscb}, we prove that our barrier function $\Gamma$ for $\kmmd$ is an LHSCB.

\subsection{Matrix monotonicity}
\label{sec:matmono:mm}

A function $f$ is \emph{matrix monotone} (or operator monotone) if $w_a \succeq w_b \succeq 0$ implies $f(w_a) \succeq f(w_b)$ for all $w_a, w_b \in \bbS^d$ for all integers $d$.
The following integral representation result is attributed to \citet{lowner1934monotone} (see e.g. \citet[Theorem 1]{kwong1989some} and \citet[Theorem L]{furuta2008concrete}).
A function $f : \bbR_{>} \to \bbR$ is matrix monotone in $\bbR_{>}$ if and only if it has the representation:
\begin{align}
f(x) 
&= \alpha + \beta x + \int_0^\infty \frac{x}{x + t} \ \diff \rho(t)
= \alpha + \beta x + \int_0^\infty \biggl( 1 - \frac{t}{x + t} \biggr) \ \diff \rho(t),
\end{align}
\label{eq:matmono:int}
where $\alpha \in \bbR, \beta \in \Rge$ and $\rho$ is a positive measure on $\bbR_{>}$ such that $\int_0^\infty (1 + t)^{-1} \ \diff \rho(t) < \infty$.

This result implies that, for a cone of squares $\cQ$ of a Jordan algebra, and for $w \in \intr(\cQ)$ with the spectral decomposition $w = \sum_{i \in \iin{d}} w_i c_i$, we have:
\begin{subequations}
\begin{align}
f(w) 
&= \sum_{i \in \iin{d}} f(w_i) c_i
\\
&= \alpha e + \beta w + \int_0^\infty \sum_{i \in \iin{d}} 
\biggl( 1 - \frac{t}{w_i + t} \biggr) c_i \ \diff \rho (t)
\\
&= \alpha e + \beta w + \int_0^\infty (e - t (w + t e)^{-1}) \ \diff \rho(t).
\end{align}
\label{eq:matmono:intQ}
\end{subequations}
This is similar to the representation in \citet[page 1520]{faybusovich2017matrix}.

\subsection{Cone definition}
\label{sec:matmono:cone}

Let $h : \bbR_> \to \bbR$ be a convex $C^3$-smooth function.
We assume that that the first derivative of $h$, $\nabla h$, is a matrix monotone function.
This also implies that $h$ is convex.
We call such functions \emph{matrix monotone derivative} (MMD) functions.

In \cref{tab:mmd}, we give some common examples of MMD functions, with abbreviated names in the first column.
We also give $\nabla h$, the domain of the convex conjugate $h^\ast$ (defined in \cref{eq:conj}), and a closed-form formula for $h^\ast$.
Due to \citet[Theorem 2.6 (L{\"o}wner-Heinz Theorem)]{carlen2010trace}, the functions $x \to \log(x)$, $x \to -x^p$ for $p \in [-1, 0]$, and $x \to x^p$ for $p \in [0, 1]$ are matrix monotone.
This implies that in \cref{tab:mmd}, each function in the $\nabla h$ column is matrix monotone.
Note that NegSqrt is equivalent to NegPower for $p = 1/2$; we highlight NegSqrt as an interesting case, for which the conjugate $h^\ast$ is a positive rescaling of the inverse function.
Note we exclude the case $p = 1$ in NegPower and Power because it is homogeneous ($h$ is linear). 
More examples of matrix monotone functions can be found in \citet{kwong1989some,furuta2008concrete}.

\begin{table}[!htb]
\centering
\begin{tabular}{ccccc}
\toprule
MMD function & $h$ & $\nabla h$ & $\dom(h^\ast)$ & $h^\ast$ \\
\midrule
NegLog & $-\log(x)$ & $-x^{-1}$ & $\bbR_>$ & $-1 - \log(x)$ \\
NegEntropy & $x \log(x)$ & $1 + \log(x)$ & $\bbR$ & $\exp(-1 - x)$ \\
NegSqrt & $-\sqrt{x}$ & $-\tfrac{1}{2} x^{-1/2}$ & $\bbR_>$ & $\tfrac{1}{4} x^{-1}$ \\
NegPower, $p \in (0,1)$ & $-x^p$ & $-p x^{p-1}$ & $\Rge$ & $-(p-1) (x/p)^q$ \\
Power, $p \in (1,2]$ & $x^p$ & $p x^{p-1}$ & $\bbR$ & $(p-1) (x_{-} / p)^q$ \\
\bottomrule
\end{tabular}

\caption{
Examples of MMD functions.
We let $q \coloneqq p / (p - 1)$, which gives $q \in (-\infty, 0)$ for $p \in (0,1)$ in the NegPower case, or $q \in [2, \infty)$ for $p \in (1,2]$ in the Power case.
We also let $x_{-} \coloneqq \max(-x, 0)$ in the Power case.
}
\label{tab:mmd}
\end{table}

\begin{figure}[!htb]
\centering
\begin{tikzpicture}
\begin{axis}[
    width=10cm,
    height=7cm,
    domain=0:2.5,
    xmin=0,
    xmax=2.5,
    ymin=-2,
    ymax=4,
    samples=400,
    xtick distance=0.5,
    ytick distance=1,
    grid=both,
    legend pos=outer north east,    
]
\addplot+[mark=none] {-ln(x)};
\addlegendentry{$-\log(x)$}
\addplot+[mark=none] {x*ln(x)};
\addlegendentry{$x \log(x)$}
\addplot+[mark=none] {-x^(0.5)};
\addlegendentry{$-\sqrt{x}$}
\addplot+[mark=none] {x^(1.5)};
\addlegendentry{$x^{3/2}$}
\end{axis}
\end{tikzpicture}
\caption{
Plots of example MMD functions.
}
\label{fig:mmd}
\end{figure}
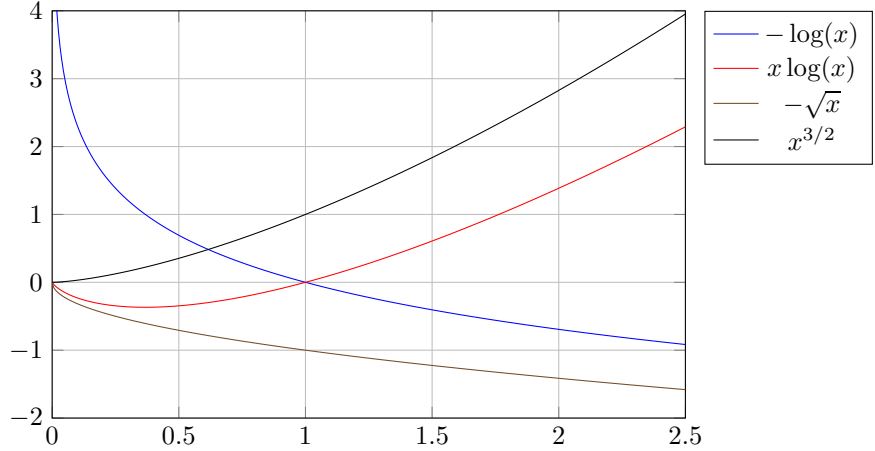

Suppose $\cQ$ is the cone of squares of a Jordan algebra $V$ with rank $d$.
Let $\varphi : \intr(\cQ) \to \bbR$ be the $C^3$-smooth function $\varphi(w) = \tr(h(w)) = \sum_{i \in \iin{d}} h(w_i)$, which is a convex separable spectral function (see \cref{sec:spectral:sep}).
As in \cref{sec:cones:per}, we let $\tu = (u, v, w) \in \cE = \bbR \times \bbR_{>} \times \intr(\cQ)$.
The function $\zeta : \cE \to \bbR$ has the form:
\begin{equation}
\zeta (\tu) \coloneqq u - v \tr(h(v^{-1} w)).
\label{eq:matmono:zeta}
\end{equation}

We define the matrix monotone derivative cone $\kmmd$, a special case of the epigraph-perspective cone $\K_p$ in \cref{eq:epiper}, as:
\begin{equation}
\kmmd \coloneqq \cl \{ 
\tu \in \cE : u \geq v \tr(h(v^{-1} w))
\},
\label{eq:matmono:K}
\end{equation}
which is a proper cone. 
Note that the negative log-determinant cone, whose barrier function we examined in \cref{sec:epiper:logdet}, is a special case of $\kmmd$ where $h$ is the NegLog function from \cref{tab:mmd}.
For the separable case, the convex conjugate $\varphi^\ast : V \to \bbR \cup \infty$ (see \cref{eq:conj}) of $\varphi$ is $\varphi^\ast(r) = \tr(h^\ast(r))$.
So from \cref{eq:epiperdual}, the dual cone is:
\begin{equation}
\kmmd^\ast 
\coloneqq \cl \{ 
\tu \in \bbR_> \times \bbR \times V : 
v \geq u \tr(h^\ast (u^{-1} w))
\}.
\label{eq:matmono:Kdu}
\end{equation}

\subsection{Derivatives of the MMD trace}
\label{sec:matmono:phid}

Suppose $w \succ 0$.
Since $\varphi(w) = \tr(h(w))$ and $\nabla h$ is matrix monotone, from \cref{eq:matmono:intQ} we can write the gradient:
\begin{equation}
\nabla \varphi(w) 
= \nabla h(w)
= \alpha e + \beta w + 
\int_0^\infty ( e - t (w + t e)^{-1} ) \ \diff \rho(t).
\label{eq:matmono:phid1}
\end{equation}
Note that $w + t e \succ 0$ for $t \geq 0$, so $(w + t e)^{-1}$ is well-defined.
Differentiating \cref{eq:matmono:phid1} in the direction $r \in V$, we have the second order directional derivative:
\begin{equation}
\nabla^2 \varphi(w) [r]
\labrel{\ref{eq:Dinv}}{=} 
\beta r +
\int_0^\infty t P((w + t e)^{-1}) r \ \diff \rho(t),
\label{eq:matmono:phid2}
\end{equation}
and the third order directional derivative:
\begin{equation}
\nabla^3 \varphi(w) [r,r]
\labrel{\ref{eq:Pinv}}{=}
-2 \int_0^\infty t P((w + t e)^{-1/2}) (P((w + t e)^{-1/2}) r)^2 \ \diff \rho(t).
\label{eq:matmono:phid3}
\end{equation}

\subsection{Self-concordant barrier}
\label{sec:matmono:lhscb}

For $\kmmd$, the LHB $\Gamma : \intr(\kmmd) \to \bbR$ from \cref{eq:lhb} has the form:
\begin{equation}
\Gamma(\tu) \coloneqq -\log(\zeta(\tu)) - \log(v) - \logdet(w).
\label{eq:matmono:bar}
\end{equation}
We describe easily-computable oracles for this $\Gamma$ in \cref{sec:epiper}, including an inverse Hessian product $\bH$ in \cref{sec:epiper:invhesssep} that is as easy to compute as the Hessian product $H$ (since $\varphi$ is a separable spectral function).

We note that \citet{faybusovich2017matrix} derive a $(1 + d)$-self-concordant barrier for the related convex (but not conic) set $\cS$:
\begin{equation}
\cS \coloneqq \cl \{ 
(u, w) \in \bbR \times \intr(\cQ) : u - \varphi(w) \geq 0 
\}.
\end{equation}
$\kmmd$ is the conic hull of $\cS$.
In \cref{lem:matmono}, we prove that our barrier $\Gamma$ in \cref{eq:matmono:bar} is self-concordant, hence it is an LHSCB for $\kmmd$ with parameter $2 + d$.
This small additive increment of one in the barrier parameter is in sharp contrast to generic conic hull results, which give barriers with a large multiplicative factor in the parameter (for example, \citet[Proposition 5.1.4]{nesterov1994interior} yields the parameter $800(1 + d)$).
Since the optimal barrier parameter for $\cl(\cE)$ is $1 + d$, our parameter cannot be reduced by more than one.

\begin{proposition}
\label{lem:matmono}
$\Gamma$ in \cref{eq:matmono:bar} is a $(2 + d)$-LHSCB for $\kmmd$ in \cref{eq:matmono:K}.
\end{proposition}

\begin{proof}
We show that $\zeta$ in \cref{eq:matmono:zeta} is $(\Rge, 1)$-compatible with the domain $\cE$, in the sense of \citet[Definition 5.1.1]{nesterov1994interior}.
This follows if (i) $\zeta$ is $C^3$-smooth on $\cE$, (ii) $\zeta$ is concave with respect to $\Rge$, (iii) for any point $\tu \in \intr (\kmmd)$ and direction $\tp = (p, q, r) \in \bbR \times \bbR \times V$ satisfying $v \pm q \geq 0$ and $w \pm r \succeq 0$ it holds that:
\begin{equation}
\nabla^3 \zeta(\tu) [\tp,\tp,\tp] 
\leq -3 \nabla^2 \zeta(\tu) [\tp,\tp].
\label{eq:matmono:compat}
\end{equation}

Suppose $v \pm q \geq 0$ and $w \pm r \succeq 0$.
As in \cref{eq:mu}, we let $\mu \coloneqq \mu(v, w) = v^{-1} w \succ 0$ and $\xi \coloneqq v^{-1} (r - q \mu)$.
From \cref{eq:zetadds}, the second and third order directional derivatives of $\zeta$ at $\tu$ in direction $\tp$ are:
\begin{subequations}
\begin{align}
\nabla^2 \zeta(\tu) [\tp,\tp]
&= -v \nabla^2 \varphi(\mu) [\xi,\xi],
\label{eq:matmono:chiderivs:D2}
\\
\nabla^3 \zeta(\tu) [\tp,\tp,\tp] 
&= -v \nabla^3 \varphi(\mu) [\xi,\xi,\xi] + 3 q \nabla^2 \varphi(\mu) [\xi,\xi].
\label{eq:matmono:chiderivs:D3}
\end{align}
\label{eq:matmono:chiderivs}
\end{subequations}
Since $\varphi$ is convex and $C^3$-smooth on $\intr(\cQ)$ by assumption, \cref{eq:matmono:zeta,eq:matmono:chiderivs} imply that $\zeta$ is concave and $C^3$-smooth on $\cE$.
It remains to show that \cref{eq:matmono:compat} holds.

For $t \geq 0$, let:
\begin{subequations}
\begin{align}
a(t) 
&\coloneqq \mu + t e \succ 0,
\\
\ba(t) 
&\coloneqq a(t)^{-1/2} \succ 0,
\\
\hxi(t) 
&\coloneqq P(\ba(t)) \xi.
\end{align}
\end{subequations}
By the integral representation result from \cref{sec:matmono:mm} \citep{lowner1934monotone}, there exists a positive measure $\rho$ and $\beta \geq 0$ such that the directional derivatives of $\varphi$ are:
\begin{subequations}
\begin{align}
\nabla^2 \varphi(\mu) [\xi,\xi] 
&\labrel{\ref{eq:matmono:phid2}}{=} 
\beta \tr(\xi^2) + \int_0^{\infty} t \tr(\hxi(t)^2) \ \diff \rho(t),
\\
\nabla^3 \varphi(\mu) [\xi,\xi,\xi] 
&\labrel{\ref{eq:matmono:phid3}}{=} 
-2 \int_0^{\infty} t \tr(\hxi(t)^3) \ \diff \rho(t).
\end{align}
\label{eq:matmono:phiderivs}
\end{subequations}
From \cref{eq:matmono:chiderivs,eq:matmono:phiderivs}, the compatibility condition \cref{eq:matmono:compat} is equivalent to nonnegativity of:
\begin{subequations}
\begin{align}
& -3 \nabla^2 \zeta(\tu) [\tp,\tp] - \nabla^3 \zeta(\tu) [\tp,\tp,\tp]
\\
&= 3 (v - q) \beta \tr(\xi^2) +
\int_0^{\infty} t ( 3 (v - q) \tr(\hxi(t)^2) - 2 v \tr(\hxi(t)^3) ) \ \diff \rho(t).
\label{eq:matmono:d2d3}
\end{align}
\end{subequations}

Since $v \geq q$, the first term in \cref{eq:matmono:d2d3} is nonnegative.
The second term (the integral) is nonnegative if for all $t \geq 0$, the following inner term is nonnegative:
\begin{equation}
3 (v - q) \tr(\hxi(t)^2) - 2 v \tr(\hxi(t)^3) 
= \langle \hxi(t)^2, 3 (v - q) e - 2 v \hxi(t) \rangle.
\label{eq:matmono:trineq}
\end{equation}
By self-duality of $\cQ$, \cref{eq:matmono:trineq} is nonnegative if $3 (v - q) e - 2 v \hxi(t) \succeq 0$, which we now prove.
For $t \geq 0$, let:
\begin{equation}
b(t) 
\coloneqq (1 - v^{-1} q) a(t) - \xi
= v^{-1} ( w - r ) + t (1 - v^{-1} q) e.
\label{eq:bt}
\end{equation}
Since $w \succeq r$ and $1 - v^{-1} q \geq 0$, we have $b(t) \succeq 0$.
Hence we have:
\begin{equation}
(1 - v^{-1} q) e - \hxi(t)
\labrel{\ref{eq:bt}}{=} 
P(\ba(t)) b(t)
\succeq 0,
\label{eq:sandwich}
\end{equation}
since $\ba(t) \succ 0$ implies $P(\ba(t))$ is an automorphism on $\cQ$ (see \citet[Page 48]{faraut1998analysis}).
Therefore:
\begin{equation}
3 (v - q) e - 2 v \hxi(t)
\labrel{\ref{eq:sandwich}}{\succeq} 
3 (v - q) e - 2 v (1 - v^{-1} q) e
= (v - q) e
\succeq 0.
\end{equation}
So \cref{eq:matmono:trineq} is nonnegative, which implies the integral term in \cref{eq:matmono:d2d3} is nonnegative.

Thus \cref{eq:matmono:d2d3} is nonnegative, so \cref{eq:matmono:compat} holds and compatibility is proved.
Now by \citet[Proposition 5.1.7]{nesterov1994interior}, $\Gamma$ is an LHSCB for $\kmmd$ with parameter $2 + d$.

\end{proof}

\section{Root-determinant cones}
\label{sec:rtdet}

In \cref{sec:rtdet:def}, we define the root-determinant cone $\krtdet$, which is the hypograph of the homogeneous nonseparable spectral root-determinant function.
After expressing the derivatives of this function in \cref{sec:rtdet:phid}, we prove that our barrier function $\Gamma$ for $\krtdet$ is an LHSCB in \cref{sec:rtdet:bar}, and we derive easily-computable barrier oracles in \cref{sec:rtdet:oracles}.

\subsection{Cone definition}
\label{sec:rtdet:def}

Suppose $\cQ$ is a cone of squares of a Jordan algebra $V$ with rank $d$.
Let $\varphi : \cQ \to \Rge$ denote the root-determinant function (or the geometric mean of the eigenvalues):
\begin{equation}
\varphi(w) \coloneqq \det(w)^{1/d} 
= \prod_{i \in \iin{d}} w_i^{1/d},
\label{eq:rtdet:phi}
\end{equation}
which is a concave, homogeneous nonseparable spectral function (see \cref{sec:spectral:gen}).
We let $\tu \coloneqq (u, w) \in \cE = \bbR \times \cQ$.
The function $\zeta : \cE \to \bbR$ has the form:
\begin{align}
\zeta(\tu) &\coloneqq \varphi(w) - u.
\label{eq:rtdet:zeta}
\end{align}

We define the root-determinant cone $\krtdet$ and its dual cone as:
\begin{subequations}
\begin{align}
\krtdet &\coloneqq 
\{ (u, w) \in \bbR \times \cQ : 
u \leq \det(w)^{1/d} \},
\label{eq:rtdet:K}
\\
\krtdet^\ast &\coloneqq 
\{ (u, w) \in \bbR_{\leq} \times \cQ : 
-d^{-1} u \leq \det(w)^{1/d} \}.
\label{eq:rtdet:Kdu}
\end{align}
\end{subequations}
We note that $\krtdet$ is a hypograph modification of the epigraph cone $\K_h$ in \cref{eq:epi}, and it is a primitive proper cone.
$\krtdet^\ast$ can be derived by modifying the steps we use to derive $\K_h^\ast$ in \cref{eq:epidualhomog} and using the convex conjugate of the negative root-determinant function.

\subsection{Derivatives of root-determinant}
\label{sec:rtdet:phid}

Suppose $w \succ 0$.
Since $\varphi(w) = \exp(d^{-1} \logdet(w))$, applying the chain rule and using \cref{eq:Dpsi} gives us the gradient:
\begin{equation}
\nabla \varphi(w)
= d^{-1} \varphi(w) w^{-1}.
\label{eq:rtdet:phid1}
\end{equation}
Let $r \in V$ and $\hr \coloneqq P(w^{-1/2}) r \in V$.
Using the product rule on \cref{eq:rtdet:phid1}, we have the second order directional derivative:
\begin{subequations}
\begin{align}
\nabla^2 \varphi(w) [r]
&= d^{-2} \langle w^{-1}, r \rangle \nabla \varphi(w) + 
d^{-1} \varphi(w) \nabla_w(w^{-1}) [r]
\\
&\labrel{\ref{eq:Dinv}}{=} 
d^{-1} \varphi(w) \tr(\hr) w^{-1} - d^{-1} \varphi(w) P(w^{-1}) r
\\
&= d^{-1} \varphi(w) ( d^{-1} \tr(\hr) w^{-1} - P(w^{-1}) r )
\label{eq:rtdet:phid2:a}
\\
&= d^{-1} \varphi(w) P(w^{-1/2}) (d^{-1} \tr(\hr) e - \hr).
\end{align}
\label{eq:rtdet:phid2}
\end{subequations}
Finally, using the product rule on \cref{eq:rtdet:phid2:a}, we have the third order directional derivative:
\begin{subequations}
\begin{align}
\begin{split}
\nabla^3 \varphi(w) [r, r]
&= d^{-1} \langle d^{-1} \tr(\hr) w^{-1} - P(w^{-1}) r, r \rangle \nabla \varphi(w) +
d^{-1} \varphi(w) ( 
\\
&\pheq d^{-1} ( 
\langle w^{-1}, r \rangle \nabla_w (\langle w^{-1}, r \rangle) + \tr(\hr) \nabla_w(w^{-1}) [r]
) - \nabla_w (P(w^{-1}) r) [r] )
\end{split}
\\
\begin{split}
&\labrel{\ref{eq:Pinv}}{=} d^{-2} \varphi(w) ( d^{-1} \tr(\hr)^2 - \tr(\hr^2) ) w^{-1} +
d^{-1} \varphi(w) (
\\
&\pheq -2 d^{-1} \tr(\hr) P(w^{-1}) r + 2 P(w^{-1/2}) (P(w^{-1/2}) r)^2 )
\end{split}
\\
&= d^{-1} \varphi(w) P(w^{-1/2}) (
d^{-1} ( d^{-1} \tr(\hr)^2 - \tr(\hr^2) ) e -
2 d^{-1} \tr(\hr) \hr + 
2 \hr^2).
\end{align}
\label{eq:rtdet:phid3}
\end{subequations}

\subsection{Self-concordant barrier}
\label{sec:rtdet:bar}

For $\krtdet$, the LHB $\Gamma : \intr(\krtdet) \to \bbR$ from \cref{eq:lhb} has the form:
\begin{equation}
\Gamma (\tu) \coloneqq 
-\log ( \zeta(\tu) ) - \logdet(w).
\label{eq:rtdet:bar}
\end{equation}
In \cref{lem:rtdet} we show that $\Gamma$ is self-concordant with parameter $1 + d$.
Since the optimal barrier parameter for $\cE$ is $d$, our parameter cannot be reduced by more than one.

\begin{proposition}
\label{lem:rtdet}
$\Gamma$ in \cref{eq:rtdet:bar} is a $(1 + d)$-LHSCB for $\krtdet$ in \cref{eq:rtdet:K}.
\end{proposition}

\begin{proof}
Note $\Psi(\tu) \coloneqq -\logdet(w)$ is a $d$-LHSCB for $\cE$.
We show that $\zeta$ in \cref{eq:rtdet:zeta} is $(\Rge, 1)$-compatible with the barrier $\Psi$ in the sense of \citet[Definition 5.1.2]{nesterov1994interior}.
Compatibility follows if (i) $\zeta$ is $C^3$-smooth on $\intr (\cE)$, (ii) concave with respect to $\Rge$, (iii) for any point $\tu \in \intr(\krtdet)$ and direction $\tp = (p, r) \in \bbR \times V$ it holds that:
\begin{equation}
\nabla^3 \zeta(\tu) [\tp, \tp, \tp] 
\leq - 3 ( \nabla^2 \Psi(\tu) [\tp, \tp] )^{1/2} \nabla^2 \zeta(\tu) [\tp, \tp] .
\label{eq:rtdet:compat}
\end{equation}

Suppose $\tu \in \intr(\krtdet)$.
From \cref{eq:rtdet:zeta}, we have:
\begin{subequations}
\begin{align}
\nabla^2 \zeta(\tu) [\tp, \tp] 
&= \nabla^2 \varphi(w) [r, r],
\\
\nabla^3 \zeta(\tu) [\tp, \tp, \tp] 
&= \nabla^3 \varphi(w) [r, r, r].
\end{align}
\label{eq:rtdet:zetadd}
\end{subequations}
Since $\varphi$ is concave and $C^3$-smooth on $\intr(\cQ)$, \cref{eq:rtdet:zetadd} implies $\zeta$ is concave and $C^3$-smooth on $\intr(\cE)$.
It remains to show that \cref{eq:rtdet:compat} holds.

Let $\sigma \in \bbR^d$ be the eigenvalues of $\hr \coloneqq P(w^{-1/2}) r$.
Then:
\begin{equation}
( \nabla^2 \Psi(\tu) [\tp, \tp] )^{1/2}
\labrel{\ref{eq:psid2}}{=} 
\tr(\hr^2)^{1/2} = \lVert \sigma \rVert.
\label{eq:rtdet:psi}
\end{equation}
Let $m_k \coloneqq d^{-1} \tr(\hr^k), \forall k \in \iin{3}$, and let $\delta_i \coloneqq \sigma_i - m_1, \forall i \in \iin{d}$.
By the formulae for variance and skewness, we have:
\begin{subequations}
\begin{align}
m_2 - m_1^2 &= 
d^{-1} \sum_{i \in \iin{d}} \delta_i^2,
\label{eq:rtdet:var}
\\
m_3 - 3 m_1 m_2 + 2 m_1^3 &=
d^{-1} \sum_{i \in \iin{d}} \delta_i^3.
\end{align}
\label{eq:rtdet:varskew}
\end{subequations}
For convenience, let $\varphi \coloneqq \varphi(w) > 0$ be a constant.
We have:
\begin{subequations}
\begin{align}
\nabla^2 \varphi(w) [r, r]
&\labrel{\ref{eq:rtdet:phid2}}{=} 
d^{-1} \varphi \langle P(w^{-1/2}) (d^{-1} \tr(\hr) e - \hr), r \rangle
\\
&= -\varphi (d^{-1} \tr(\hr^2) - d^{-2} \tr(\hr)^2)
\\
&= -\varphi (m_2 - m_1^2)
\\
&\labrel{\ref{eq:rtdet:var}}{=} 
-d^{-1} \varphi \sum_{i \in \iin{d}} \delta_i^2
\leq 0.
\end{align}
\label{eq:rtdet:phidd2}
\end{subequations}
Similarly:
\begin{subequations}
\begin{align}
\nabla^3 \varphi(w) [r, r, r]
&\labrel{\ref{eq:rtdet:phid3}}{=} 
d^{-1} \varphi \langle
d^{-1} ( d^{-1} \tr(\hr)^2 - \tr(\hr^2) ) e - 2 d^{-1} \tr(\hr) \hr + 2 \hr^2
, \hr \rangle
\\
&= d^{-1} \varphi (
d^{-1} ( d^{-1} \tr(\hr)^2 - \tr(\hr^2) ) \tr(\hr) - 2 d^{-1} \tr(\hr) \tr(\hr^2) + 2 \tr(\hr^3)
)
\\
&= \varphi ( m_1^3 - 3 m_1 m_2 + 2 m_3 )
\\
&= \varphi ( 3 m_1 ( m_2 - m_1^2 ) + 2 ( m_3 - 3 m_1 m_2 + 2 m_1^3 ) )
\\
&\labrel{\ref{eq:rtdet:varskew}}{=} 
d^{-1} \varphi \sum_{i \in \iin{d}} ( 3 m_1 \delta_i^2 + 2 \delta_i^3 )
\\
&= d^{-1} \varphi \sum_{i \in \iin{d}} \delta_i^2 ( m_1 + 2 \sigma_i ).
\end{align}
\label{eq:rtdet:phidd3}
\end{subequations}

Finally, using \cref{eq:rtdet:psi,eq:rtdet:phidd2,eq:rtdet:phidd3} the compatibility condition \cref{eq:rtdet:compat} is equivalent to nonnegativity of:
\begin{subequations}
\begin{align}
& -\nabla^3 \zeta(\tu) [\tp, \tp, \tp] - 
3 ( \nabla^2 \Psi(\tu) [\tp, \tp] )^{1/2} \nabla^2 \zeta(\tu) [\tp, \tp]
\\
&= -d^{-1} \varphi \sum_{i \in \iin{d}} \delta_i^2 ( m_1 + 2 \sigma_i ) +
3 \lVert \sigma \rVert d^{-1} \varphi \sum_{i \in \iin{d}} \delta_i^2
\\
&= d^{-1} \varphi \sum_{i \in \iin{d}} \delta_i^2 
( \lVert \sigma \rVert - m_1 + 2 (\lVert \sigma \rVert - \sigma_i) ).
\label{eq:rtdet:sciii}
\end{align}
\end{subequations}
Clearly, $d^{-1} \varphi \delta_i^2 \geq 0$ and $\sigma_i \leq \lVert \sigma \rVert$ for all $i \in \iin{d}$.
We have $m_1 \leq d^{-1} \lVert \sigma \rVert_1 \leq d^{-1/2} \lVert \sigma \rVert \leq \lVert \sigma \rVert$.
Hence \cref{eq:rtdet:sciii} is nonnegative.

Thus \cref{eq:rtdet:compat} holds and compatibility is proved.
Now by \citet[Proposition 5.1.7]{nesterov1994interior}, $\Gamma$ is a $(1 + d)$-LHSCB for $\krtdet$.

\end{proof}

\subsection{Evaluating barrier oracles}
\label{sec:rtdet:oracles}

Using the derivatives of $\varphi$ from \cref{sec:rtdet:phid}, we derive easily-computable oracles for the LHSCB \cref{eq:rtdet:bar}.
Let $\tu \in \intr(\krtdet)$ and $\tp = (p, r) \in \bbR \times V$.
For convenience, let $\varphi \coloneqq \varphi(w) > 0$ be a constant. 
We define the scalar constants:
\begin{subequations}
\begin{align}
\eta &\coloneqq d^{-1} \varphi \zeta^{-1},
\\
\theta &\coloneqq 1 + \eta,
\\
\chi &\coloneqq -\zeta^{-1} p + \eta \tr(\hr),
\\
\tau &\coloneqq \chi - d^{-1} \tr(\hr),
\\
\upsilon &\coloneqq \tr(\hr^2) - d^{-1} \tr(\hr)^2.
\end{align}
\end{subequations}
Note that:
\begin{subequations}
\begin{align}
\nabla_u (\zeta^{-1}) &= \zeta^{-2},
\\
\nabla_u \eta &= \zeta^{-1} \eta,
\\
\nabla_u \chi &= \zeta^{-1} \chi,
\\
\nabla_w (\zeta^{-1}) 
&= -\zeta^{-2} \nabla \varphi(w) 
\labrel{\ref{eq:rtdet:phid1}}{=}
-\zeta^{-1} \eta w^{-1},
\\
\nabla_w \eta 
&= \eta (d^{-1} - \eta) w^{-1},
\\
\nabla_w \chi
&\labrel{\ref{eq:Dinv}}{=}
\zeta^{-1} \eta p w^{-1} + \eta (d^{-1} - \eta) \tr(\hr) w^{-1} - \eta P(w^{-1}) r
\\
&= -\eta (\tau w^{-1} + P(w^{-1}) r).
\end{align}
\end{subequations}

The gradient of $\Gamma$ in \cref{eq:rtdet:bar} is:
\begin{subequations}
\begin{align}
g_u
&= \zeta^{-1},
\\
g_w
&= -\zeta^{-1} \nabla \varphi(w) - w^{-1}
\\
&= -\theta w^{-1}.
\end{align}
\label{eq:rtdet:g}
\end{subequations}
Differentiating \cref{eq:rtdet:g}, the Hessian product is:
\begin{subequations}
\begin{align}
H_u 
&= \nabla_u g_u p + \nabla_u g_w [r]
\\
&= \zeta^{-2} p - \zeta^{-1} \eta \tr(\hr)
\\
&= -\zeta^{-1} \chi,
\\
H_w 
&= \nabla_w g_u p + \nabla_w g_w [r]
\\
&\labrel{\ref{eq:Dinv}}{=}
-\zeta^{-1} \eta p w^{-1} - \tr(\hr) \eta (d^{-1} - \eta) w^{-1} + \theta P(w^{-1}) r 
\\
&= P(w^{-1/2}) (\eta \tau e + \theta \hr).
\end{align}
\label{eq:rtdet:H}
\end{subequations}
Differentiating \cref{eq:rtdet:H}, the the third order directional derivative is:
\begin{subequations}
\begin{align}
T_u 
&= \nabla_u H_u p + \nabla_u H_w [r]
\\
&= -2 \zeta^{-2} p \chi +
\zeta^{-1} \eta (\tau \tr(\hr) + \tr(\hr^2)) + \zeta^{-1} \eta \tr(\hr) \chi
\\
&= \zeta^{-1} ( 2 \chi^2 + \eta \upsilon ),
\\
T_w 
&= \nabla_w H_u p + \nabla_w H_w [r]
\\
\begin{split}
&\labrel{\ref{eq:Pinv}}{=}
\zeta^{-1} \eta p \chi w^{-1} + \zeta^{-1} \eta p (\tau w^{-1} + P(w^{-1}) r) - 
\eta \tau P(w^{-1}) r + {}
\\
&\pheq \tr(\hr) \tau \eta (d^{-1} - \eta) w^{-1} + 
\eta ( -\eta (\tau \tr(\hr) + \tr(\hr^2)) + d^{-1} \tr(\hr^2)) w^{-1} + {}
\\
&\pheq \eta (d^{-1} - \eta) \tr(\hr) P(w^{-1}) r
- 2 \theta P(w^{-1/2}) \hr^2
\end{split}
\\
\begin{split}
&= \eta ( 
-\chi \tau + \zeta^{-1} p \chi + (d^{-1} - \eta) (\tr(\hr) \tau + \tr(\hr^2))
) w^{-1} + {}
\\
&\pheq \eta ( -\tau + \zeta^{-1} p + (d^{-1} - \eta) \tr(\hr) ) P(w^{-1}) r 
- 2 \theta P(w^{-1/2}) \hr^2
\end{split}
\\
&= P(w^{-1/2}) (
\eta (-2 \chi \tau + (d^{-1} - \eta) \upsilon) e -
2 \eta \tau \hr - 2 \theta \hr^2 ).
\end{align}
\end{subequations}

In \cref{lem:rtdet:invhess} below, we give a closed-form inverse Hessian product operator.
This operator \cref{eq:rtdet:invhess} a similar structure to the Hessian product operator \cref{eq:rtdet:H}, except that it applies $P(w^{1/2})$ instead of $P(w^{-1/2})$.

\begin{lemma}
\label{lem:rtdet:invhess}

Letting $\chr \coloneqq P(w^{1/2}) r \in V$, the inverse Hessian product is:
\begin{subequations}
\begin{align}
\bH_u 
&= (\zeta^2 + d^{-1} \varphi^2) p + d^{-1} \varphi \tr(\chr),
\\
\bH_w 
&= P(w^{1/2}) (
d^{-1} (\varphi p + \eta \theta^{-1} \tr(\chr)) e +
\theta^{-1} \chr
).
\end{align}
\label{eq:rtdet:invhess}
\end{subequations}
\end{lemma}

\begin{proof}
Note that the Hessian operator \cref{eq:rtdet:H} is a positive definite linear operator, so it has a unique inverse linear operator.
We show that $(\nabla^2 \Gamma)^{-1} (\nabla^2 \Gamma [\tp]) = \tp$.
Into \cref{eq:rtdet:invhess}, we substitute the values from \cref{eq:rtdet:H} i.e. $p = H_u = -\zeta^{-1} \chi$ and $r = H_w = P(w^{-1/2}) (\eta \tau e + \theta \hr)$.
Since $P(w^{1/2}) = P(w^{-1/2})^{-1}$, we have:
\begin{subequations}
\begin{align}
\chr 
&= P(w^{1/2}) H_w 
= \eta \tau e + \theta \hr,
\\
\tr(\chr) 
&= d \eta \tau + \theta \tr(\hr).
\end{align}
\end{subequations}
We have:
\begin{subequations}
\begin{align}
\bH_u 
&= (\zeta^2 + d^{-1} \varphi^2) (-\zeta^{-1} \chi) + d^{-1} \varphi (d \eta \tau + \theta \tr(\hr))
\\
&= -\zeta \chi + \varphi \eta (\tau - \chi) + d^{-1} \varphi \theta \tr(\hr)
\\
&= -\zeta (\chi - \eta \tr(\hr))
\\
&= p,
\end{align}
\end{subequations}
and:
\begin{subequations}
\begin{align}
\bH_w 
&= P(w^{1/2}) (
d^{-1} (-\varphi \zeta^{-1} \chi + \eta \theta^{-1} (d \eta \tau + \theta \tr(\hr)) ) e +
\theta^{-1} (\eta \tau e + \theta \hr) )
\\
&= P(w^{1/2}) (
\eta (-\tau + \eta \theta^{-1} \tau + \theta^{-1} \tau) e +
\hr )
\\
&= P(w^{1/2}) (\hr)
\\
&= r.
\end{align}
\end{subequations}
Hence \cref{eq:rtdet:invhess} is the unique inverse operator of \cref{eq:rtdet:H}.

\end{proof}

We note the polynomial-like structure of the oracles.
In particular, the $w$ components of the $g$, $H$, and $T$ oracles are computed by applying $P(w^{-1/2})$ to a polynomial in $\hr$, of degree zero for $g$, degree one for $H$, and degree two for $T$.
Analogously to $H$, its inverse $\bH$ is computed by applying $P(w^{1/2})$ to a polynomial of degree one in $\chr$.
This structure leads to simple, efficient, and numerically-stable implementations.
We also note the structural similarity (ignoring constants) between the $u$ and $w$ components of these oracles and those of the negative log-determinant barrier in \cref{sec:epiper:logdet}.
In both cases, the oracles can be computed without an explicit eigendecomposition if it is possible to apply $P(w^{1/2})$ and $P(w^{-1/2})$ directly.
For example for $V = \bbS^d$ and $V = \bbH^d$, only a Cholesky factorization of $w$ is needed.

\section{Examples and computational testing}
\label{sec:testing}

We outline our implementations of the MMD cone and the log-determinant and root-determinant cones in Hypatia in \cref{sec:testing:hyp}.
In \crefrange{sec:testing:nonparam}{sec:testing:classquant}, we present example problems with simple, natural formulations (NFs) in terms of these cones.
Using techniques we describe in \cref{sec:testing:form}, we construct equivalent extended formulations (EFs) that can be recognized by MOSEK 9 or ECOS. 
Our computational benchmarks follow the methodology we describe in \cref{sec:testing:meth} and show that Hypatia often solves the NFs much more efficiently than Hypatia, MOSEK, or ECOS solve the EFs.
Finally, in \cref{sec:testing:invhessprod}, we exemplify the computational impact of efficient oracle procedures by comparing the performance of our closed-form inverse Hessian product formula in \cref{eq:invhessprodsep} with that of a naive direct solve using the explicit Hessian matrix.

\subsection{Hypatia solver}
\label{sec:testing:hyp}

Hypatia's generic cone interface allows specifying a \emph{vectorized} proper cone $\K \subset \bbR^q$ for some dimension $q$.
The nonnegative cone is already in vectorized format.
For the real symmetric PSD cone $\Sge^d$, we use the standard \emph{svec} transformation, which rescales and stores only the elements of the matrix triangle in a vector of dimension $d (d + 1) / 2$. 
For the complex Hermitian PSD cone $\Hge^d$, we perform a modified \emph{svec} transformation to a $d^2$-dimensional vector, storing each real diagonal element as a single element and each complex off-diagonal element in the triangle as two (rescaled) consecutive real elements (the real part followed by the imaginary part).
These transformations preserve inner products and the self-duality of cones of squares.

We adapt these transformations to enable vectorization of spectral cones. 
For example, for the epigraph-perspective cone $\K_p$ in \cref{eq:epiper}, the vectorization is $(u, v, \vect(w)) \in \bbR^{2 + q}$, where $\vect(w) \in \bbR^q$ is the appropriate vectorization of $w \in \cQ$.
Fortunately, the dual cone of this vectorized cone is the analogous vectorization of the dual cone $\K_p^\ast$ in \cref{eq:epiperdual}.

Hypatia's primal general conic form over a vector of variables $x$ minimizes a linear function $\langle c, x \rangle$ subject to linear equality constraints $b - A x = 0$ and conic constrains $h - G x \in \K$.
Here, the vectorized proper cone $\K$ is specified as a Cartesian product $\K_1 \times \cdots \times \K_K$ of $K$ proper cones recognized by Hypatia.
This means for each $\K_i$, we must have either $\K_i$ or $\K_i^\ast$ defined through Hypatia's cone interface.

For the domains $\bbR^d$, $\bbS^d$, and $\bbH^d$, we implement vectorizations of the MMD cone $\kmmd$, the log-determinant cone $\klogdet$, and the root-determinant cone ($\krtdet$).\footnote{
Our $\klogdet$ implementation is for the hypograph of the log-determinant, rather than the epigraph of negative log-determinant considered in \cref{sec:epiper:logdet}.
This only requires minor changes to the oracle derivations and the LHSCB proof from \cref{sec:matmono:lhscb} for validity.
}
This allows the user to model with these cones or their dual cones in Hypatia.
As we discuss at the end of \cref{sec:epiper:logdet,sec:rtdet:oracles}, for $\klogdet$ and $\krtdet$ the oracle procedures are quite specialized, for example we compute a Cholesky factorization rather than an eigendecomposition for the $\bbS^d$ and $\bbH^d$ domains.

For $\kmmd$, we predefine the MMD functions in \cref{tab:mmd} (e.g. \emph{NegEntropy}).
Recall that $\kmmd^\ast$ in \cref{eq:matmono:Kdu} is defined using the convex conjugate of the MMD function; in the examples below we suffix the MMD function names with \emph{Conj} (e.g. \emph{NegEntropyConj}) to indicate use of the convex conjugate function and $\kmmd^\ast$.
We write \emph{NegLogdet} or \emph{NegRtdet} for the negative log-determinant or negative root-determinant function, the epigraph of which we represent using $\klogdet$ or $\krtdet$.
In our examples, we choose not to use $\krtdet^\ast$ or $\klogdet^\ast$ (or equivalently, $\kmmd^\ast$ with \emph{NegLogConj}), because these particular dual cones provide little additional modeling power over their primal cones.

\subsection{Natural and extended formulations}
\label{sec:testing:form}

To assess the computational value of our new cones and efficient oracles, we compare the performance of Hypatia on NFs over $\kmmd$, $\klogdet$, and $\krtdet$ against that of other conic IPM solvers on equivalent EFs.
ECOS \citep{domahidi2013ecos} is another open-source conic IPM solver, but it only supports nonnegative, second-order, and three-dimensional exponential cones.
MOSEK version 9 \citep{mosek2020modeling} is a commercial conic IPM solver that supports the same cones as ECOS as well as three-dimensional power cones and real symmetric PSD cones.
We call the cones supported by MOSEK 9 the \emph{standard cones}.
Hypatia currently supports around two dozen cone types (not counting dual cones) \citep{coey2021performance}, including $\kmmd$, $\klogdet$, $\krtdet$, and (generalizations of) the standard cones.
To build the standard cone EFs, we use a variety of formulation techniques, some of which we discuss and analyze in \citet{coey2021solving}.\footnote{
The EFs build automatically via the functions in \url{https://github.com/chriscoey/Hypatia.jl/blob/master/examples/spectral_functions_JuMP.jl}.}

For $V = \bbR^d$, our EFs are constructed as follows.
The EFs for \emph{NegLog}, \emph{NegEntropy}, \emph{NegEntropyConj}, and \emph{NegRtdet} use $d$ exponential cones.
The EFs for \emph{NegSqrt} and \emph{NegSqrtConj} use $d$ three-dimensional second-order cones.
The EFs for \emph{Power}, \emph{NegPower}, \emph{PowerConj}, and \emph{NegPowerConj} use $d$ power cones.
The example in \cref{sec:testing:nonparam} uses $V = \bbR^d$.

For $V = \bbS^d$, our EFs are constructed as follows.
For most spectral functions, we adapt the EF from \citet[Proposition 4.2.1]{ben2001lectures}, which requires using an EF from the $V = \bbR^d$ case for the corresponding spectral function, and adding constraints on the sum of the $i$ largest eigenvalues of a matrix for each $i \in \iin{d}$.
This is a large formulation with many additional variables and PSD constraints.
For \emph{NegLog} and \emph{NegRtdet}, we use a much simpler EF from \citet[Example 18.d]{ben2001lectures}.
Since \emph{NegSqrtConj} is a scaling of the inverse function (see \cref{tab:mmd}), a Schur complement representation allows us to use an EF with one PSD cone constraint.
For $V = \bbH^d$, we reformulate any complex PSD cone constraint to a real PSD cone constraint with twice the side dimension \citep[Section 6.2.7]{mosek2020modeling}.
The examples in \cref{sec:testing:expdesign,sec:testing:centpoly} use $V = \bbS^d$ and the example in \cref{sec:testing:classquant} uses $V = \bbH^d$.

\subsection{Computational methodology}
\label{sec:testing:meth}

We perform all instance generation, computational experiments, and results analysis with Ubuntu 21.10, Julia 1.8.0-DEV.862, and Hypatia 0.5.3.\footnote{
Benchmark scripts and instructions for reproducing and analyzing results are available at \url{https://github.com/chriscoey/Hypatia.jl/tree/master/benchmarks/natvsext}. 
A raw output CSV file and detailed results tables are at \url{https://github.com/chriscoey/Hypatia.jl/wiki}.
}
We use dedicated hardware with an AMD Ryzen 9 3950X 16-core processor (32 threads) and 128GB of RAM.
For each example problem in \crefrange{sec:testing:nonparam}{sec:testing:classquant}, we generate random instances of a range of sizes, using JuMP 0.21.10 and MathOptInterface v0.9.22.
All instances are primal-dual feasible, so we expect solvers to return optimality certificates.

We use the conic PDIPM solvers in MOSEK version 9 and ECOS version 2.0.5 (with no features disabled).
Hypatia uses a particular default algorithmic implementation that we describe in \citet{coey2021performance} (the combined directions method with the QR-Cholesky linear system procedure).
We limit each solver to 16 threads and set a solve time limit of 1800 seconds.
We set relative feasibility and optimality gap tolerances to $10^{-7}$ and absolute optimality gap tolerances to $10^{-10}$.

For each instance, the relative difference between the objective values of the formulation/solver combinations that converge never exceeds $10^{-4}$.
For each instance/formulation/solver combination that returns a solution, we measure the maximum violation $\epsilon$ of the primal-dual optimality conditions in \citet[Equation 23]{coey2021solving}.
In \crefrange{fig:nonparam}{fig:classquant}, we plot solve times in seconds against an instance size parameter, excluding solves for which $\epsilon > 10^{-5}$.
Hypatia-NF (i.e. Hypatia solving the NF) is faster than any EF solver (Hypatia-EF, MOSEK-EF, ECOS-EF) across all instance sizes and spectral functions tested for each example, and always scales to larger sizes.

\subsection{Examples and results}
\label{sec:testing:ex}

\subsubsection{Nonparametric distribution estimation}
\label{sec:testing:nonparam}

Suppose we have a random variable $X$ taking values in the finite set $\{\alpha_i\}_{i \in \iin{d}}$.
We seek a probability distribution $\rho \in \bbR^d$ that minimizes a convex spectral function $\varphi$, given some prior information expressed with $d/2$ linear equality constraints.
Adapting \citet[Section 7.2]{boyd2004convex}, the problem is:
\begin{subequations}
\begin{align}
\textstyle\min_{\rho \in \bbR^d} \quad \varphi(\rho) &:
\\
\tr(\rho) & = d,
\\
A \rho & = b.
\end{align}
\label{eq:nonparam}
\end{subequations}
For four spectral functions $\varphi$ on $\Rge^d$ (with EFs that ECOS can recognize) and a range of sizes $d$, we build random instances of \cref{eq:nonparam}.
The solver timings are summarized in \cref{fig:nonparam}.
Note that for \emph{NegRtdet}, no solve times are plotted for MOSEK-EF because the optimality condition violations $\epsilon$ are too large (see \cref{sec:testing:meth}); tightening MOSEK's tolerance options improves these violations, though in either case MOSEK-EF is significantly slower than Hypatia-NF.
We do not plot results for \emph{NegLogdet} (the $\klogdet$ formulation using the specialized oracles from \cref{sec:epiper:logdet}) as they are nearly identical to the results for $\kmmd$/\emph{NegLog}; however, the efficiency benefits of \emph{NegLogdet} are realized for the matrix domain in \cref{sec:testing:expdesign}.

\begin{figure}[!htb]
\centering
\begin{tikzpicture}
\pgfplotstableread[col sep = comma]{csvs/nonparametricdistr_plot_NegRtdet.csv} \csvnegrtdet
\pgfplotstableread[col sep = comma]{csvs/nonparametricdistr_plot_NegLogSSF.csv} \csvneglog
\pgfplotstableread[col sep = comma]{csvs/nonparametricdistr_plot_NegEntropySSF.csv} \csvnegentropy
\pgfplotstableread[col sep = comma]{csvs/nonparametricdistr_plot_NegSqrtSSF.csv} \csvnegsqrt
\begin{groupplot}[
    group style={group size=2 by 2, horizontal sep=1.2cm, vertical sep=1.3cm, ylabels at=edge left},
    width=0.415\textwidth,
    height = 6cm,
    scale only axis,
    ymax = 2000,
    xmin = 0,
    enlargelimits=false,
    ymode=log,
    xmajorgrids=true,
    ymajorgrids=true,
    ylabel=time (s),
    ylabel shift=-5pt,
    xlabel=$d$,
    x label style={at={(axis description cs:1.0,0)}, anchor=north, yshift=-1pt},
    title style={yshift=-5pt},
    xtick = {0,10000,20000},
    scaled x ticks = false,
    every tick label/.append style={font=\small},
    every x tick label/.append style={yshift=-1pt},
    every axis plot/.append style={
        line width=0.5pt,
        mark options={solid},
        mark size=2pt,
        },
    ]
\nextgroupplot[
    title = NegRtdet, 
    ]
\addplot [blue, mark=o] table [x=d, y=nat_Hypatia] from \csvnegrtdet;
\addplot [black, mark=x] table [x=d, y=vecext_Hypatia] from \csvnegrtdet;
\addplot [olive, mark=square] table [x=d, y=vecext_ECOS] from \csvnegrtdet;
\nextgroupplot[
    title = NegLog, 
    ]
\addplot [blue, mark=o] table [x=d, y=nat_Hypatia] from \csvneglog;
\addplot [black, mark=x] table [x=d, y=vecext_Hypatia] from \csvneglog;
\addplot [red, mark=triangle] table [x=d, y=vecext_Mosek] from \csvneglog;
\addplot [olive, mark=square] table [x=d, y=vecext_ECOS] from \csvneglog;
\nextgroupplot[
    title = NegEntropy, 
    ]
\addplot [blue, mark=o] table [x=d, y=nat_Hypatia] from \csvnegentropy;
\addplot [black, mark=x] table [x=d, y=vecext_Hypatia] from \csvnegentropy;
\addplot [red, mark=triangle] table [x=d, y=vecext_Mosek] from \csvnegentropy;
\addplot [olive, mark=square] table [x=d, y=vecext_ECOS] from \csvnegentropy;
\nextgroupplot[
    title = NegSqrt, 
    legend style={
        legend columns=-1,
        fill=none,
        draw=black,
        /tikz/every even column/.append style={column sep=0.5cm},
        /tikz/every odd column/.append style={column sep=0.1cm},
        },
    legend to name=leg3,
    ]
\addplot [blue, mark=o] table [x=d, y=nat_Hypatia] from \csvnegsqrt;
\addplot [black, mark=x] table [x=d, y=vecext_Hypatia] from \csvnegsqrt;
\addplot [red, mark=triangle] table [x=d, y=vecext_Mosek] from \csvnegsqrt;
\addplot [olive, mark=square] table [x=d, y=vecext_ECOS] from \csvnegsqrt;
\legend{Hypatia-NF,Hypatia-EF,MOSEK-EF,ECOS-EF};
\end{groupplot}
\path (group c1r2.south east) -- node[below, yshift=-0.7cm]{\pgfplotslegendfromname{leg3}} (group c2r2.south west);
\end{tikzpicture}
\caption{
\nameref{sec:testing:nonparam} solver performance.
}
\label{fig:nonparam}
\end{figure}

\subsubsection{Experiment design}
\label{sec:testing:expdesign}

We formulate a continuous relaxation of the experiment design problem, similar to \citet[Section 7.5]{boyd2004convex}.
The variable $\rho \in \bbR^{2d}$ is the number of trials to run for each of $2d$ experiments that are useful for estimating a vector in $\bbR^d$.
The experiments are described by the columns of $V \in \bbR^{d \times 2d}$ and we require that $2d$ experiments are performed.
We minimize a convex spectral function of the information matrix:
\begin{subequations}
\begin{align}
\textstyle\min_{\rho \in \bbR^{2d}} \quad \varphi(V \Diag(\rho) V') &:
\\
\tr(\rho) & = 2d,
\\
\rho & \geq 0,
\end{align}
\label{eq:expdesign}
\end{subequations}
where $V'$ is the transpose of $V$ and $\Diag(\rho)$ is the diagonal matrix of $\rho$.
For four different $\varphi$ on $\Sge^d$ and various $d$, we build random instances of \cref{eq:expdesign}.
The solver timings are summarized in \cref{fig:expdesign}.
Since ECOS does not support $\Sge^d$, we only compare with MOSEK.
The \emph{Hypatia-NegLogdet} curve indicates that Hypatia with $\klogdet$ is somewhat more efficient than Hypatia with the equivalent $\kmmd$/\emph{NegLog} formulation; this is due to our oracle specializations in \cref{sec:epiper:logdet} and our implementation using a Cholesky factorization rather than an eigendecomposition.

\begin{figure}[!htb]
\centering
\begin{tikzpicture}
\pgfplotstableread[col sep = comma]{csvs/experimentdesign_plot_NegRtdet.csv} \csvnegrtdet
\pgfplotstableread[col sep = comma]{csvs/experimentdesign_plot_NegLogSSF.csv} \csvneglog
\pgfplotstableread[col sep = comma]{csvs/experimentdesign_plot_NegSqrtSSFConj.csv} \csvnegsqrtconj
\pgfplotstableread[col sep = comma]{csvs/experimentdesign_plot_NegPower01SSF.csv} \csvnegpower
\begin{groupplot}[
    group style={group size=2 by 2, horizontal sep=1.2cm, vertical sep=1.3cm, ylabels at=edge left},
    width=0.415\textwidth,
    height = 6cm,
    scale only axis,
    ymax = 2000,
    xmin = 0,
    enlargelimits=false,
    ymode=log,
    xmajorgrids=true,
    ymajorgrids=true,
    ylabel=time (s),
    ylabel shift=-5pt,
    xlabel=$d$,
    x label style={at={(axis description cs:1.0,0)}, anchor=north, yshift=-1pt},
    title style={yshift=-5pt},
    scaled x ticks = false,
    every tick label/.append style={font=\small},
    every x tick label/.append style={yshift=-1pt},
    every axis plot/.append style={
        line width=0.5pt,
        mark options={solid},
        mark size=2pt,
        },
    ]
\nextgroupplot[
    title = NegRtdet, 
    xtick = {0,200,...,800},
    ]
\addplot [blue, mark=o] table [x=d, y=nat_Hypatia] from \csvnegrtdet;
\addplot [black, mark=x] table [x=d, y=ext_Hypatia] from \csvnegrtdet;
\addplot [red, mark=triangle] table [x=d, y=ext_Mosek] from \csvnegrtdet;
\nextgroupplot[
    title = NegLog, 
    xtick = {0,200,...,800},
    legend style={
        legend columns=-1,
        fill=none,
        draw=black,
        /tikz/every even column/.append style={column sep=0.5cm},
        /tikz/every odd column/.append style={column sep=0.1cm},
        },
    legend to name=leg2,
    ]
\addplot [teal, mark=*] table [x=d, y=natlog_Hypatia] from \csvneglog;
\addplot [blue, mark=o] table [x=d, y=nat_Hypatia] from \csvneglog;
\addplot [black, mark=x] table [x=d, y=ext_Hypatia] from \csvneglog;
\addplot [red, mark=triangle] table [x=d, y=ext_Mosek] from \csvneglog;
\legend{Hypatia-NegLogdet,Hypatia-NF,Hypatia-EF,MOSEK-EF};
\nextgroupplot[
    title = NegSqrtConj, 
    xtick = {0,200,400,600},
    ]
\addplot [blue, mark=o] table [x=d, y=nat_Hypatia] from \csvnegsqrtconj;
\addplot [black, mark=x] table [x=d, y=ext_Hypatia] from \csvnegsqrtconj;
\addplot [red, mark=triangle] table [x=d, y=ext_Mosek] from \csvnegsqrtconj;
\nextgroupplot[
    title = NegPower($1/3$), 
    xtick = {0,200,...,800},
    ]
\addplot [blue, mark=o] table [x=d, y=nat_Hypatia] from \csvnegpower;
\addplot [black, mark=x] table [x=d, y=ext_Hypatia] from \csvnegpower;
\addplot [red, mark=triangle] table [x=d, y=ext_Mosek] from \csvnegpower;
\end{groupplot}
\path (group c1r2.south east) -- node[below, yshift=-0.7cm]{\pgfplotslegendfromname{leg2}} (group c2r2.south west);
\end{tikzpicture}
\caption{
\nameref{sec:testing:expdesign} solver performance.
}
\label{fig:expdesign}
\end{figure}

\subsubsection{Central polynomial Gram matrix}
\label{sec:testing:centpoly}

Suppose we have a polynomial of degree $2k$ in $m$ variables.
Let $L = \binom{m+k}{m}$ and $U = \binom{m+2k}{m}$, and let $b \in \bbR^U$ be the monomial coefficients of the polynomial.
We seek a Gram matrix $\rho \in \bbS^L$ corresponding to $b$ \citep[Lemma 3.33]{parrilo2012chapter} that minimizes a convex spectral function $\varphi$:
\begin{subequations}
\begin{align}
\textstyle\min_{\rho \in \bbS^L} \quad \varphi(\rho) & :
\\
C \vect(\rho) & = b,
\end{align}
\label{eq:centpoly}
\end{subequations}
where the matrix $C \in \bbR^{U \times L(L+1)/2}$ maps the Gram matrix to the (lower-dimensional) polynomial coefficient space.
We build random instances of \cref{eq:centpoly}, varying $m \in \{1,4\}$ and $k$ (depending on $m$).
Recall from \cref{tab:mmd} that \emph{ConjNegEntr} and \emph{ConjPower-1.5} are defined on $\bbS^d$, but \emph{NegEntr} and \emph{MatPower12(1.5)} are only defined on $\Sge^d$, which implicitly requires that $b$ be a sum of squares polynomial and hence globally nonnegative.
The solver timings are summarized in \cref{fig:centpoly} (a log-log plot). 

\begin{figure}[!htb]
\centering
\begin{tikzpicture}
\pgfplotstableread[col sep = comma]{csvs/centralpolymat_plot_NegEntropySSF_1.csv} \csvnegentropya
\pgfplotstableread[col sep = comma]{csvs/centralpolymat_plot_NegEntropySSF_4.csv} \csvnegentropyb
\pgfplotstableread[col sep = comma]{csvs/centralpolymat_plot_NegEntropySSFConj_1.csv} \csvnegentropyconja
\pgfplotstableread[col sep = comma]{csvs/centralpolymat_plot_NegEntropySSFConj_4.csv} \csvnegentropyconjb
\pgfplotstableread[col sep = comma]{csvs/centralpolymat_plot_Power12SSF_1.csv} \csvpowera
\pgfplotstableread[col sep = comma]{csvs/centralpolymat_plot_Power12SSF_4.csv} \csvpowerb
\pgfplotstableread[col sep = comma]{csvs/centralpolymat_plot_Power12SSFConj_1.csv} \csvpowerconja
\pgfplotstableread[col sep = comma]{csvs/centralpolymat_plot_Power12SSFConj_4.csv} \csvpowerconjb
\begin{groupplot}[
    group style={group size=2 by 2, horizontal sep=1.2cm, vertical sep=1.3cm, ylabels at=edge left},
    width=0.415\textwidth,
    height = 6cm,
    scale only axis,
    ymin = 0.07,
    ymax = 2000,
    xmin = 1,
    xmax = 200,
    enlargelimits=false,
    xmode=log, 
    ymode=log,
    xmajorgrids=true,
    ymajorgrids=true,
    ylabel=time (s),
    ylabel shift=-5pt,
    xlabel=$k$,
    x label style={at={(axis description cs:1.0,0)}, anchor=north, yshift=-1pt},
    title style={yshift=-5pt},
    every tick label/.append style={font=\small},
    every x tick label/.append style={yshift=-1pt},
    xtick={1,10,100},
    xticklabels={1,10,100},
    every axis plot/.append style={
        line width=0.5pt,
        mark options={solid},
        mark size=2pt,
        },
    ]
\nextgroupplot[
    title = NegEntropy, 
    ]
\addplot [blue, mark=o] table [x=k, y=nat_Hypatia] from \csvnegentropya;
\addplot [blue, mark=o, dashed] table [x=k, y=nat_Hypatia] from \csvnegentropyb;
\addplot [black, mark=x] table [x=k, y=ext_Hypatia] from \csvnegentropya;
\addplot [black, mark=x, dashed] table [x=k, y=ext_Hypatia] from \csvnegentropyb;
\addplot [red, mark=triangle] table [x=k, y=ext_Mosek] from \csvnegentropya;
\addplot [red, mark=triangle, dashed] table [x=k, y=ext_Mosek] from \csvnegentropyb;
\nextgroupplot[
    title = NegEntropyConj, 
    ]
\addplot [blue, mark=o] table [x=k, y=nat_Hypatia] from \csvnegentropyconja;
\addplot [blue, mark=o, dashed] table [x=k, y=nat_Hypatia] from \csvnegentropyconjb;
\addplot [black, mark=x] table [x=k, y=ext_Hypatia] from \csvnegentropyconja;
\addplot [black, mark=x, dashed] table [x=k, y=ext_Hypatia] from \csvnegentropyconjb;
\addplot [red, mark=triangle] table [x=k, y=ext_Mosek] from \csvnegentropyconja;
\addplot [red, mark=triangle, dashed] table [x=k, y=ext_Mosek] from \csvnegentropyconjb;
\nextgroupplot[
    title = Power12($1.5$), 
    ]
\addplot [blue, mark=o] table [x=k, y=nat_Hypatia] from \csvpowera;
\addplot [blue, mark=o, dashed] table [x=k, y=nat_Hypatia] from \csvpowerb;
\addplot [black, mark=x] table [x=k, y=ext_Hypatia] from \csvpowera;
\addplot [black, mark=x, dashed] table [x=k, y=ext_Hypatia] from \csvpowerb;
\addplot [red, mark=triangle] table [x=k, y=ext_Mosek] from \csvpowera;
\addplot [red, mark=triangle, dashed] table [x=k, y=ext_Mosek] from \csvpowerb;
\nextgroupplot[
    title = Power12Conj($1.5$), 
    legend style={
        legend columns=-1,
        fill=none,
        draw=black,
        /tikz/every even column/.append style={column sep=0.5cm},
        /tikz/every odd column/.append style={column sep=0.1cm},
        },
    legend to name=leg1,
    ]
\addlegendimage{black, line legend}
\addlegendentry{$m = 1$}
\addlegendimage{black, dashed, line legend}
\addlegendentry{$m = 4$}
\addlegendimage{blue, mark=o, only marks}
\addlegendentry{Hypatia-NF}
\addlegendimage{black, mark=x, only marks}
\addlegendentry{Hypatia-EF}
\addlegendimage{red, mark=triangle, only marks}
\addlegendentry{MOSEK-EF}
\addplot [blue, mark=o] table [x=k, y=nat_Hypatia] from \csvpowerconja;
\addplot [blue, mark=o, dashed] table [x=k, y=nat_Hypatia] from \csvpowerconjb;
\addplot [black, mark=x] table [x=k, y=ext_Hypatia] from \csvpowerconja;
\addplot [black, mark=x, dashed] table [x=k, y=ext_Hypatia] from \csvpowerconjb;
\addplot [red, mark=triangle] table [x=k, y=ext_Mosek] from \csvpowerconja;
\addplot [red, mark=triangle, dashed] table [x=k, y=ext_Mosek] from \csvpowerconjb;
\end{groupplot}
\path (group c1r2.south east) -- node[below, yshift=-0.7cm]{\pgfplotslegendfromname{leg1}} (group c2r2.south west);
\end{tikzpicture}
\caption{
\nameref{sec:testing:centpoly} solver performance.
}
\label{fig:centpoly}
\end{figure}

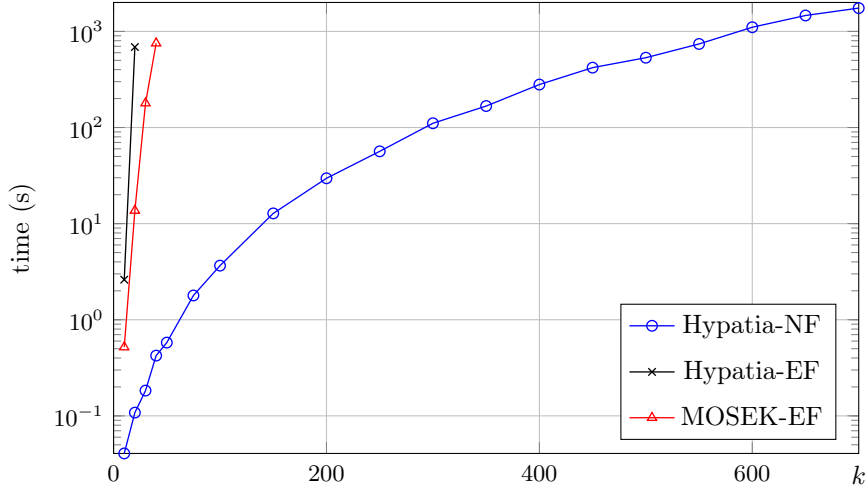
\subsubsection{Classical-quantum channel capacity}
\label{sec:testing:classquant}

We compute the capacity of a classical-quantum channel, adapting the formulation from \citet[Example 2.16]{sutter2015efficient} and \citet[Section 3.1]{fawzi2018efficient}.
The variable $\rho \in \bbR^d$ is a probability distribution on the $d$-dimensional input alphabet.
For $i \in \iin{d}$, let $P_i \in \Hge^d$ be fixed density matrices satisfying $\tr(P_i) = 1$.
Letting $\varphi$ represent the trace of \emph{NegEntropy} on $\Hge^d$, the formulation is:
\begin{subequations}
\begin{align}
\textstyle\min_{\rho \in \bbR^d} \quad 
\varphi \bigl( \tsum{i \in \iin{d}} \rho_i P_i \bigr) - \sum_{i \in \iin{d}} \rho_i \varphi(P_i) &:
\\
\tr(\rho) & = 1,
\\
\rho & \geq 0.
\end{align}
\label{eq:classquant}
\end{subequations}
We generate random instances of \cref{eq:classquant}, varying $d$.
The solver timings are summarized in \cref{fig:classquant}. 

\begin{figure}[!htb]
\centering
\begin{tikzpicture}
\pgfplotstableread[col sep = comma]{csvs/classicalquantum_plot.csv} \csvclassquant
\begin{axis}[
    width=0.6\textwidth,
    height = 6cm,
    scale only axis,
    ymax = 2000,
    ymode=log,
    xmin = 0,
    xmajorgrids=true,
    ymajorgrids=true,
    enlargelimits=false,
    ylabel=time (s),
    xlabel=$k$,
    xtick={0,200,400,600},
    x label style={at={(axis description cs:1.0,0)}, anchor=north, yshift=-1pt},
    every tick label/.append style={font=\small},
    every x tick label/.append style={yshift=-1pt},
    every axis plot/.append style={
        line width=0.5pt,
        mark options={solid},
        mark size=2pt,
        },
    legend pos=south east,
    legend style={row sep=5pt},
    ]
\addplot [blue, mark=o] table [x=d, y=nat_Hypatia] from \csvclassquant;
\addplot [black, mark=x] table [x=d, y=ext_Hypatia] from \csvclassquant;
\addplot [red, mark=triangle] table [x=d, y=ext_Mosek] from \csvclassquant;
\legend{Hypatia-NF,Hypatia-EF,MOSEK-EF};
\end{axis}
\end{tikzpicture}
\caption{
\nameref{sec:testing:classquant} solver performance.
}
\label{fig:classquant}
\end{figure}
\subsection{Inverse Hessian product oracle}
\label{sec:testing:invhessprod}

To illustrate the importance of our efficient and numerically stable oracle procedures, we compare the performance of two different approaches to computing the inverse Hessian product oracle $\bH$ in \cref{eq:oracles:bh} for $\kmmd$ cones.
The naive approach is to compute the explicit Hessian matrix, perform a Cholesky factorization, and use a direct linear solve.
Alternatively, we derive a closed-form formula for $\bH$ in \cref{eq:invhessprodsep}, since $\kmmd$ is a special case of $\K_p$ with a separable spectral function.
This formula is essentially as easy to compute as the Hessian product oracle $H$ in \cref{eq:hessprod} (which does not use an explicit Hessian matrix).
In \cref{tab:complexity}, we compare the worst-case memory and time complexities for these procedures.

\begin{table}[!htb]
\centering
\begin{tabular}{clllll}
\toprule
& & \multicolumn{2}{c}{closed-form formula} & \multicolumn{2}{c}{factorize and solve} \\
\cmidrule(lr){3-4} \cmidrule(lr){5-6}
$V$ & $\dim(\kmmd)$ & memory & time & memory & time \\
\midrule
$\bbR^d$ & $O(d)$ &
$O(d)$ & $O(d)$ &
$O(d^2)$ & $O(d^3)$
\\
$\bbS^d$ or $\bbH^d$ & $O(d^2)$ &
$O(d^2)$ & $O(d^3)$ & 
$O(d^4)$ & $O(d^6)$
\\
\bottomrule
\end{tabular}
\caption{
Cone dimension and worst-case complexities for the two inverse Hessian product procedures.
}
\label{tab:complexity}
\end{table}

To compare the practical performance of these procedures, we perform computational experiments using Hypatia.
We first solve NF instances of a range of sizes for the examples from \cref{sec:testing:nonparam} (with $V = \bbR^d$) and \cref{sec:testing:expdesign} (with $V = \bbS^d$), using $\kmmd$ with the \emph{NegEntropy} function.
For each instance, at Hypatia's final PDIPM iterate, we take the direction $r = g$ (i.e. the gradient oracle in \cref{eq:oracles:g} at the iterate) and compute $\bH$ for this direction using both procedures.
To measure the numerical accuracy of each procedure, we compute $\epsilon \coloneqq \lvert 1 - \nu^{-1} \langle \bH, g \rangle \rvert$, which measures the violation of a particular identity \citep[Equation 2.5]{nesterov1997self} satisfied by a logarithmically homogeneous function such as the LHSCB $\Gamma$ with parameter $\nu = 2 + d$.
We also time each procedure, excluding Hessian memory allocation time for the factorization-based procedure.

\begin{figure}[!htb]
\centering
\begin{tikzpicture}
\pgfplotstableread[col sep = comma]{csvs/invhess_nonparam.csv} \csvinvhessnonparam
\pgfplotstableread[col sep = comma]{csvs/invhess_expdesign.csv} \csvinvhessexpdesign
\begin{groupplot}[
    group style={group size=2 by 1, horizontal sep=1.3cm, ylabels at=edge left},
    width=0.41\textwidth,
    height = 6cm,
    scale only axis,
    ymin = 2e-12,
    ymax = 100,
    enlargelimits=false,
    xmode=log, 
    ymode=log,
    xmajorgrids=true,
    ymajorgrids=true,
    ylabel shift=-5pt,
    ylabel=time (s) or violation,
    xlabel=$d$,
    x label style={at={(axis description cs:1.0,0)}, anchor=north, yshift=-1pt},
    title style={yshift=-5pt},
    every tick label/.append style={font=\small},
    every x tick label/.append style={yshift=-1pt},
    every axis plot/.append style={
        line width=0.5pt,
        mark options={solid},
        mark size=2pt,
        },
    ]
\nextgroupplot[
    title = {\nameref{sec:testing:nonparam}: $V = \bbR^d$}, 
    xmin=3,
    xmax=10000,
    xtick={1,10,100,1000,10000},
    xticklabels={,10,100,1000,},
    ]
\addplot [blue, mark=o] table [x=d, y=t1] from \csvinvhessnonparam;
\addplot [red, mark=triangle] table [x=d, y=t2] from \csvinvhessnonparam;
\addplot [blue, mark=o, dashed] table [x=d, y=v1] from \csvinvhessnonparam;
\addplot [red, mark=triangle, dashed] table [x=d, y=v2] from \csvinvhessnonparam;
\nextgroupplot[
    title = {\nameref{sec:testing:expdesign}: $V = \bbS^d$}, 
    xmin=2,
    xmax=200,
    xtick={1,10,100,1000},
    xticklabels={,10,100,},
    legend style={
        legend columns=-1,
        fill=none,
        draw=black,
        /tikz/every even column/.append style={column sep=0.5cm},
        /tikz/every odd column/.append style={column sep=0.1cm},
        },
    legend to name=leg4,
    ]
\addlegendimage{black, line legend}
\addlegendentry{time (s)}
\addlegendimage{black, dashed, line legend}
\addlegendentry{violation}
\addlegendimage{blue, mark=o, only marks}
\addlegendentry{closed-form formula}
\addlegendimage{red, mark=triangle, only marks}
\addlegendentry{factorize and solve}
\addplot [blue, mark=o] table [x=d, y=t1] from \csvinvhessexpdesign;
\addplot [red, mark=triangle] table [x=d, y=t2] from \csvinvhessexpdesign;
\addplot [blue, mark=o, dashed] table [x=d, y=v1] from \csvinvhessexpdesign;
\addplot [red, mark=triangle, dashed] table [x=d, y=v2] from \csvinvhessexpdesign;
\end{groupplot}
\path (group c1r1.south east) -- node[below, yshift=-0.7cm]{\pgfplotslegendfromname{leg4}} (group c2r1.south west);
\end{tikzpicture}
\caption{
For instances of two examples using $\kmmd$ with \emph{NegEntropy}, the speed and logarithmic homogeneity condition violation (at the final iterate) for the two inverse Hessian product procedures.
}
\label{fig:invhessprod}
\end{figure}
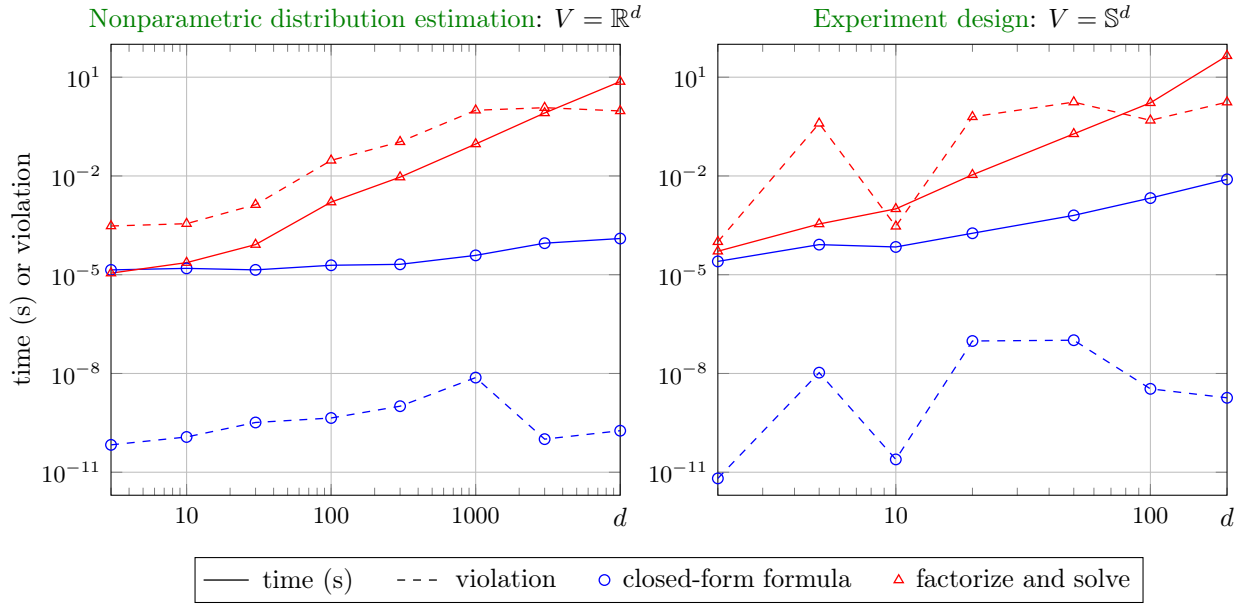

Our results are displayed in \cref{fig:invhessprod}.
The Cholesky factorization fails at $d = 3000$ for the $\bbR^d$ example and at $d = 20, 50, 200$ for the $\bbS^d$ example; when this occurs, Hypatia uses a Bunch-Kaufman factorization as a fallback (note Julia calls performant OpenBLAS routines for the Cholesky and Bunch-Kaufman factorizations).
Note that we loosen the convergence tolerances specified in \cref{sec:testing:meth} by a factor of 100, so that the factorization-based procedure fails less often.
These comparisons demonstrate that our closed-form formula generally allows computing $\bH$ faster and with greater numerical accuracy.
Also, the closed-form procedure is much more memory efficient than the factorization-based procedure, as it never forms an explicit Hessian matrix.

\bibliographystyle{plainnat}
\bibliography{refs}

\begin{thebibliography}{39}
\providecommand{\natexlab}[1]{#1}
\providecommand{\url}[1]{\texttt{#1}}
\expandafter\ifx\csname urlstyle\endcsname\relax
  \providecommand{\doi}[1]{doi: #1}\else
  \providecommand{\doi}{doi: \begingroup \urlstyle{rm}\Url}\fi

\bibitem[Andersen et~al.(2011)Andersen, Dahl, Liu, Vandenberghe, Sra, Nowozin,
  and Wright]{andersen2011interior}
Martin Andersen, Joachim Dahl, Zhang Liu, Lieven Vandenberghe, S~Sra,
  S~Nowozin, and SJ~Wright.
\newblock Interior-point methods for large-scale cone programming.
\newblock \emph{Optimization for Machine Learning}, 5583, 2011.

\bibitem[Baes(2007)]{baes2007convexity}
Michel Baes.
\newblock Convexity and differentiability properties of spectral functions and
  spectral mappings on {E}uclidean {J}ordan algebras.
\newblock \emph{Linear algebra and its applications}, 422\penalty0
  (2-3):\penalty0 664--700, 2007.

\bibitem[Ben-Tal and Nemirovski(2001)]{ben2001lectures}
Aharon Ben-Tal and Arkadi Nemirovski.
\newblock \emph{Lectures on modern convex optimization: analysis, algorithms,
  and engineering applications}.
\newblock SIAM, 2001.

\bibitem[Borchers(1999)]{borchers1999csdp}
Brian Borchers.
\newblock {CSDP}, a {C} library for semidefinite programming.
\newblock \emph{Optimization Methods and Software}, 11\penalty0 (1-4):\penalty0
  613--623, 1999.

\bibitem[Boyd et~al.(2004)Boyd, Boyd, and Vandenberghe]{boyd2004convex}
Stephen Boyd, Stephen~P Boyd, and Lieven Vandenberghe.
\newblock \emph{Convex optimization}.
\newblock Cambridge university press, 2004.

\bibitem[Carlen(2010)]{carlen2010trace}
Eric Carlen.
\newblock Trace inequalities and quantum entropy: an introductory course.
\newblock \emph{Entropy and the quantum}, 529:\penalty0 73--140, 2010.

\bibitem[Coey et~al.(2021{\natexlab{a}})Coey, Kapelevich, and
  Vielma]{coey2021performance}
Chris Coey, Lea Kapelevich, and Juan~Pablo Vielma.
\newblock Performance enhancements for a generic conic interior point
  algorithm.
\newblock \emph{arXiv preprint arXiv:2107.04262}, 2021{\natexlab{a}}.

\bibitem[Coey et~al.(2021{\natexlab{b}})Coey, Kapelevich, and
  Vielma]{coey2021solving}
Chris Coey, Lea Kapelevich, and Juan~Pablo Vielma.
\newblock Solving natural conic formulations with {H}ypatia.jl.
\newblock \emph{arXiv preprint arXiv:2005.01136}, 2021{\natexlab{b}}.

\bibitem[Davis(1957)]{davis1957all}
Chandler Davis.
\newblock All convex invariant functions of {H}ermitian matrices.
\newblock \emph{Archiv der Mathematik}, 8\penalty0 (4):\penalty0 276--278,
  1957.

\bibitem[Deng(2011)]{deng2011generalization}
Chun~Yuan Deng.
\newblock A generalization of the {S}herman--{M}orrison--{W}oodbury formula.
\newblock \emph{Applied Mathematics Letters}, 24\penalty0 (9):\penalty0
  1561--1564, 2011.

\bibitem[Domahidi et~al.(2013)Domahidi, Chu, and Boyd]{domahidi2013ecos}
Alexander Domahidi, Eric Chu, and Stephen Boyd.
\newblock {ECOS:} an {SOCP} solver for embedded systems.
\newblock In \emph{2013 European Control Conference (ECC)}, pages 3071--3076.
  IEEE, 2013.

\bibitem[Faraut and Koranyi(1998)]{faraut1998analysis}
Jacques Faraut and Adam Koranyi.
\newblock Analysis on symmetric cones.
\newblock \emph{Bull. Amer. Math. Soc}, 35:\penalty0 77--86, 1998.

\bibitem[Fawzi and Fawzi(2018)]{fawzi2018efficient}
Hamza Fawzi and Omar Fawzi.
\newblock Efficient optimization of the quantum relative entropy.
\newblock \emph{Journal of Physics A: Mathematical and Theoretical},
  51\penalty0 (15):\penalty0 154003, 2018.

\bibitem[Faybusovich and Tsuchiya(2017)]{faybusovich2017matrix}
Leonid Faybusovich and Takashi Tsuchiya.
\newblock Matrix monotonicity and self-concordance: how to handle quantum
  entropy in optimization problems.
\newblock \emph{Optimization Letters}, 11\penalty0 (8):\penalty0 1513--1526,
  2017.

\bibitem[Faybusovich and Zhou(2021)]{faybusovich2021long}
Leonid Faybusovich and Cunlu Zhou.
\newblock Long-step path-following algorithm for quantum information theory:
  Some numerical aspects and applications.
\newblock \emph{Numerical Algebra, Control \& Optimization}, 2021.

\bibitem[Furuta(2008)]{furuta2008concrete}
Takayuki Furuta.
\newblock Concrete examples of operator monotone functions obtained by an
  elementary method without appealing to {L}{\"o}wner integral representation.
\newblock \emph{Linear algebra and its applications}, 429\penalty0
  (5-6):\penalty0 972--980, 2008.

\bibitem[Grant and Boyd(2014)]{grant2014cvx}
Michael Grant and Stephen Boyd.
\newblock {CVX}: {MATLAB} software for disciplined convex programming, version
  2.1, 2014.

\bibitem[Grant et~al.(2006)Grant, Boyd, and Ye]{grant2006disciplined}
Michael Grant, Stephen Boyd, and Yinyu Ye.
\newblock Disciplined convex programming.
\newblock In \emph{Global optimization}, pages 155--210. Springer, 2006.

\bibitem[Kwong(1989)]{kwong1989some}
Man~Kam Kwong.
\newblock Some results on matrix monotone functions.
\newblock \emph{Linear Algebra and Its Applications}, 118:\penalty0 129--153,
  1989.

\bibitem[Lasserre(1998)]{lasserre1998homogeneous}
Jean~B Lasserre.
\newblock Homogeneous functions and conjugacy.
\newblock \emph{Journal of Convex Analysis}, 5:\penalty0 397--404, 1998.

\bibitem[L{\"o}wner(1934)]{lowner1934monotone}
Karl L{\"o}wner.
\newblock {\"U}ber monotone matrixfunktionen.
\newblock \emph{Mathematische Zeitschrift}, 38\penalty0 (1):\penalty0 177--216,
  1934.

\bibitem[{MOSEK ApS}(2020)]{mosek2020modeling}
{MOSEK ApS}.
\newblock {Modeling Cookbook} revision 3.2.1, 2020.
\newblock URL \url{https://docs.mosek.com/modeling-cookbook/index.html}.

\bibitem[Nesterov and Nemirovskii(1994)]{nesterov1994interior}
Y.~Nesterov and A.~Nemirovskii.
\newblock \emph{Interior-point polynomial algorithms in convex programming}.
\newblock Studies in Applied Mathematics. Society for Industrial and Applied
  Mathematics, 1994.

\bibitem[Nesterov and Todd(1997)]{nesterov1997self}
Yu~E Nesterov and Michael~J Todd.
\newblock Self-scaled barriers and interior-point methods for convex
  programming.
\newblock \emph{Mathematics of Operations research}, 22\penalty0 (1):\penalty0
  1--42, 1997.

\bibitem[Nesterov et~al.(1996)Nesterov, Todd, and Ye]{nesterov1996infeasible}
Yu~E Nesterov, Michael~J Todd, and Yinyu Ye.
\newblock Infeasible-start primal-dual methods and infeasibility detectors for
  nonlinear programming problems.
\newblock Technical report, Cornell University Operations Research and
  Industrial Engineering, 1996.

\bibitem[Nesterov(2012)]{nesterov2012towards}
Yuri Nesterov.
\newblock Towards non-symmetric conic optimization.
\newblock \emph{Optimization Methods and Software}, 27\penalty0 (4-5):\penalty0
  893--917, 2012.

\bibitem[Papp and Alizadeh(2013)]{papp2013semidefinite}
David Papp and Farid Alizadeh.
\newblock Semidefinite characterization of sum-of-squares cones in algebras.
\newblock \emph{SIAM Journal on Optimization}, 23\penalty0 (3):\penalty0
  1398--1423, 2013.

\bibitem[Parrilo(2012)]{parrilo2012chapter}
Pablo~A Parrilo.
\newblock Chapter 3: Polynomial optimization, sums of squares, and
  applications.
\newblock In \emph{Semidefinite Optimization and Convex Algebraic Geometry},
  pages 47--157. SIAM, 2012.

\bibitem[Permenter et~al.(2017)Permenter, Friberg, and
  Andersen]{permenter2017solving}
Frank Permenter, Henrik~A Friberg, and Erling~D Andersen.
\newblock Solving conic optimization problems via self-dual embedding and
  facial reduction: a unified approach.
\newblock \emph{SIAM Journal on Optimization}, 27\penalty0 (3):\penalty0
  1257--1282, 2017.

\bibitem[Rockafellar(2015)]{rockafellar2015convex}
Ralph~Tyrell Rockafellar.
\newblock \emph{Convex analysis}.
\newblock Princeton university press, 2015.

\bibitem[Sendov(2007)]{sendov2007higher}
Hristo~S Sendov.
\newblock The higher-order derivatives of spectral functions.
\newblock \emph{Linear algebra and its applications}, 424\penalty0
  (1):\penalty0 240--281, 2007.

\bibitem[Serrano(2015)]{serrano2015algorithms}
Santiago~Akle Serrano.
\newblock \emph{Algorithms for unsymmetric cone optimization and an
  implementation for problems with the exponential cone}.
\newblock PhD thesis, Stanford University, 2015.

\bibitem[Skajaa and Ye(2015)]{skajaa2015homogeneous}
Anders Skajaa and Yinyu Ye.
\newblock A homogeneous interior-point algorithm for nonsymmetric convex conic
  optimization.
\newblock \emph{Mathematical Programming}, 150\penalty0 (2):\penalty0 391--422,
  2015.

\bibitem[Sun and Sun(2008)]{sun2008lowner}
Defeng Sun and Jie Sun.
\newblock L{\"o}wner's operator and spectral functions in {E}uclidean {J}ordan
  algebras.
\newblock \emph{Mathematics of Operations Research}, 33\penalty0 (2):\penalty0
  421--445, 2008.

\bibitem[Sutter et~al.(2015)Sutter, Sutter, Esfahani, and
  Renner]{sutter2015efficient}
David Sutter, Tobias Sutter, Peyman~Mohajerin Esfahani, and Renato Renner.
\newblock Efficient approximation of quantum channel capacities.
\newblock \emph{IEEE Transactions on Information Theory}, 62\penalty0
  (1):\penalty0 578--598, 2015.

\bibitem[Vieira(2007)]{vieira2007jordan}
Manuel~V.C. Vieira.
\newblock \emph{Jordan Algebraic approach to symmetric optimization}.
\newblock PhD thesis, Faculty of Electrical Engineering, Mathematics and
  Computer Science, TU Delft, NL 2628 CD, Delft, The Netherlands, November
  2007.

\bibitem[Vieira(2016)]{vieira2016derivatives}
Manuel~VC Vieira.
\newblock Derivatives of eigenvalues and {J}ordan frames.
\newblock \emph{Numerical Algebra, Control \& Optimization}, 6\penalty0
  (2):\penalty0 115, 2016.

\bibitem[Yamashita et~al.(2003)Yamashita, Fujisawa, and
  Kojima]{yamashita2003implementation}
Makoto Yamashita, Katsuki Fujisawa, and Masakazu Kojima.
\newblock Implementation and evaluation of {SDPA} 6.0 (semidefinite programming
  algorithm 6.0).
\newblock \emph{Optimization Methods and Software}, 18\penalty0 (4):\penalty0
  491--505, 2003.

\bibitem[Zhang(2004)]{zhang2004new}
Shuzhong Zhang.
\newblock A new self-dual embedding method for convex programming.
\newblock \emph{Journal of Global Optimization}, 29\penalty0 (4):\penalty0
  479--496, 2004.

\end{thebibliography}


\end{document}